\documentstyle{article}

\input amssym.def
\input amssym

\newif\ifComments \newif\ifTwelve

\Commentstrue

\DeclareFontShape{OT1}{cmr}{m}{sc}
 {<5><6><7><8>cmcsc8 <9>cmcsc9%
 <10><10.95><12><14.4><17.28><20.74><24.88> cmcsc10}{}

 \let\tenline\newline
\ifTwelve  \let\tenline\relax \fi

\ifx\mathrm\undefined\let\mathrm\bf\fi
\ifx\mathbf\undefined\let\mathbf\bf\fi

\let\le\leqslant \let\leq\leqslant
\let\ge\geqslant \let\geq\geqslant
\let\llee\preccurlyeq

\def\hybrid{\topmargin 0pt \oddsidemargin 0pt \headheight 0pt \headsep 0pt
 \textwidth 160truemm           
 \textheight 9.4truein          
 \marginparwidth 0.0in
 \parskip 0pt plus 1pt \jot = 1.5ex}

\catcode`\@=11
\def\marginnote#1{}

\newcount\hour
\newcount\minute
\newtoks\amorpm
\hour=\time\divide\hour by60
\minute=\time{\multiply\hour by60 \global\advance\minute by-\hour}
\edef\standardtime{{\ifnum\hour<12 \global\amorpm={am}\else
 \global\amorpm={pm}\advance\hour by-12\fi \ifnum\hour=0 \hour=12 \fi
 \number\hour:\ifnum\minute<10 0\fi\number\minute\the\amorpm}}
\edef\militarytime{\number\hour:\ifnum\minute<10 0\fi\number\minute}

\def\draftlabel#1{{\@bsphack\if@filesw{\let\thepage\relax
 \xdef\@gtempa{\write\@auxout{\string
 \newlabel{#1}{{\@currentlabel}{\thepage}}}}}\@gtempa
 \if@nobreak\ifvmode\nobreak\fi\fi\fi\@esphack}\gdef\@eqnlabel{#1}}
\def\@eqnlabel{}
\def\@vacuum{}
\def\draftmarginnote#1{\marginpar{\raggedright\scriptsize\tt#1}}

\def\draft{\oddsidemargin -.5truein
 \def\@oddfoot{\sl preliminary draft\hfil\rm\thepage\hfil\sl\today\quad
 \militarytime} \let\@evenfoot\@oddfoot \overfullrule 3pt
 \let\label=\draftlabel \let\marginnote=\draftmarginnote
 \def\@eqnnum{(\theequation)\rlap{\kern\marginparsep\tt\@eqnlabel}%
 \global\let\@eqnlabel\@vacuum}}

\makeatletter
\newdimen\normalarrayskip               
\newdimen\minarrayskip                  
\normalarrayskip\baselineskip
\minarrayskip\jot
\newif\ifold \oldtrue 
\def\arraymode{\ifold\relax\else\displaystyle\fi} 
\def\eqnumphantom{\phantom{(\theequation)}}     
\def\@arrayskip{\ifold\baselineskip\z@\lineskip\z@\else
 \baselineskip\minarrayskip\lineskip2\minarrayskip\fi}
\def\@arrayclassz{\ifcase \@lastchclass \@acolampacol \or
\@ampacol \or \or \or \@addamp \or \@acolampacol \or \@firstampfalse \@acol\fi
\edef\@preamble{\@preamble\ifcase\@chnum\hfil$\relax\arraymode\@sharp$\hfil
 \or$\relax\arraymode\@sharp$\hfil\or \hfil$\relax\arraymode\@sharp$\fi}}
\def\@array[#1]#2{\setbox\@arstrutbox=\hbox{\vrule
 height\arraystretch \ht\strutbox depth\arraystretch \dp\strutbox
 width\z@}\@mkpream{#2}\edef\@preamble{\halign \noexpand\@halignto
\bgroup \tabskip\z@ \@arstrut \@preamble \tabskip\z@ \cr}%
\let\@startpbox\@@startpbox \let\@endpbox\@@endpbox
 \if #1t\vtop \else \if#1b\vbox \else \vcenter \fi\fi \bgroup\let\par\relax
 \let\@sharp##\let\protect\relax\@arrayskip\@preamble}
\def\eqnarray{\stepcounter{equation}\let\@currentlabel=\theequation
 \global\@eqnswtrue\global\@eqcnt\z@\tabskip\@centering\let\\ \@eqncr
 $$\halign to \displaywidth\bgroup\eqnumphantom\@eqnsel\hskip\@centering
 $\displaystyle\tabskip\z@{##}$&\global\@eqcnt\@ne \hskip 2\arraycolsep
 $\displaystyle\arraymode{##}$\hfil&\global\@eqcnt\tw@ \hskip 2\arraycolsep
 $\displaystyle\tabskip\z@{##}$\hfil\tabskip\@centering&{##}\tabskip\z@\cr}
\makeatother

\def\bea{\begin{eqnarray}}
\def\eea{\end{eqnarray}}
\def\beq{\begin{equation}}
\def\eeq{\end{equation}}
\def\nn{\nonumber}
\def\ot{\otimes}
\def\sk#1{\left({#1}\right)}

\def\tr{\mathop{\rm tr}\nolimits}

\def\Re{\mathop{\rm Re\>}\nolimits}
\def\Im{\mathop{\rm Im\>}\nolimits}
\def\Asym{\mathop{\rm Asym\>}\nolimits}
\def\sgn{\mathop{\rm sgn\>}\nolimits}
\def\RR{\Bbb{R}}
\def\ZZ{\Bbb{Z}}
\def\CC{\Bbb{C}}
\def\r#1{{\rm(\ref{#1})}}
\def\Ga#1{\Gm\left(#1\right)}
\def\fract#1#2{{\displaystyle{#1\over #2}}}
\def\stackreb#1#2{\mathop{#1\>}\limits_{#2}}
\def\res#1{\stackreb{\rm res}{#1}}
\def\lim#1{\stackreb{\rm lim}{#1}}
\def\Res#1{\stackreb{\rm Res}{#1}}
\let\dis=\displaystyle
\def\ee{{\rm e}}
\renewcommand{\theequation}{{\thesection}.{\arabic{equation}}}

\def\H{{\cal H}}
\let\z y
\let\b z
\let\g p
\let\v t
\let\u u

\def\Fc{{\cal F}}
\def\F{\widehat\Fc}
\def\Fq{\F_{\!q}}
\def\Fcq{\Fc_{\!q}}
\def\Fll{\Fc^{(\ell)}}
\def\Fol{\Fc^{\ot\ell}}
\def\Fl{\Fc[\ell\)]}
\def\Fcl{\Fc_{\!cl}[\ell\)]}
\def\Fw{\textstyle{\bigwedge^\ell}\Fc}
\def\Fww{\textstyle{\bigwedge^\ell}\Fll}
\let\h\hbar
\let\scr=\scriptstyle
\def\Aya{\widehat{{\cal D}Y}_\h({\frak{sl}}_2)}

\makeatletter

\def\Cup{\bigcup\limits} \def\Sum{\sum\limits} 
\let\Int\int \def\int{\Int\limits}
\def\sltwo{\frak{sl}_2}
\def\refskips{\parskip 0pt plus .1pt\partopsep0pt\frenchspacing}
\def\vsk#1>{\vskip#1\baselineskip} \def\Par{\endgraf\vsk.5>} \let\nt\noindent
\def\vv#1>{\vadjust{\vsk#1>}\ignorespaces}
\let\vp\vphantom \let\hp\hphantom \def\qqquad{\qquad\quad}
\def\lsym#1{#1\allowbreak\ldots\relax#1\allowbreak} \def\lc{\lsym,}
\def\ftext#1{{\let\thefootnote\relax\footnotetext{\vsk-.8>\noindent #1}}}
\def\cd{\ \cdot\ }

\def\;{\relax\ifmmode\mskip\thickmuskip\else\kern.27777em\fi}
\def\!{\relax\ifmmode\mskip-\thinmuskip\else\negthinspace\fi}
\def\>{{\;\!}} \def\){{\;\;\!\!\!}} \def\]{{\;\!\!}}

\def\gad#1{\global\advance#1\@ne}

\def\textindent#1{\indent\llap{#1\enspace}\ignorespaces}

\def\textitem#1{\par\hangindent1.5\parindent
 \hglue-.5\parindent\textindent{\hp{\upshape b)}\llap{#1}}}

\newcount\itemlet
\def\newbi{\itemlet 96} \newbi
\def\bitem{\gad\itemlet \par\hangindent1.5\parindent
 \hglue-.5\parindent\textindent{\upshape\rlap{\char\the\itemlet}\hp{b})}}
\def\atem{\newbi\bitem}

\newcount\itemrm
\def\newitrm{\itemrm 0}
\def\iitem{\gad\itemrm \par\hangindent1.5\parindent\hglue-.5\parindent
 \textindent{\upshape\hp{v}\llap{\romannumeral\the\itemrm})}}
\def\Item{\newitrm\iitem}

\let\gm\gamma \let\Gm\Gamma
\let\dl\delta \let\Dl\Delta
\let\ep\varepsilon
\let\ka\kappa
\let\la\lambda

\def\Cc{{\cal C}}
\def\Ch{\rlap{$\)\widehat{\phantom{\Cc}}\]$}\Cc}
\def\Ec{{\cal E}}
\def\Gc{{\cal G}} \def\Gcl{\Gc_{cl}}
\def\Qc{{\cal Q}}

\def\Mbar{{\;\smash{\overline{\!\]\]M\!}\,}\vp M}}
\let\tsize\textstyle \let\dsize\displaystyle
\let\der\partial
\def\@oddfoot{\footnotesize\hfil\rm\thepage\hfil}
\let\@evenfoot\@oddfoot

\def\qed{\hbox{}\nobreak\hfill\nobreak{\m@th$\,\square$}}
\def\Ucup{{\tsize\bigcup}}

\def\Sb{{\bf S}}

\def\lg{\frak l}
\def\mg{\frak m}
\def\Xg{\frak X}

\def\@xthm#1#2{\@begintheorem{#2}{\csname the#1\endcsname.}\ignorespaces}
\def\@ythm#1#2[#3]{\@opargbegintheorem{#2}{\csname the#1\endcsname}{#3.}%
 \ignorespaces}
\def\@begintheorem#1#2{\trivlist
 \item[\hskip \labelsep{\bfseries #1\ #2}]\slshape}

\def\s@ct@n{\setcounter{prop}{0}\setcounter{equation}{0}%
 \@startsection{section}{2}{\z@}{-3.25ex\@plus -1ex \@minus -.2ex}%
 {1.5ex \@plus .2ex}{\normalfont\large\bfseries}}
\def\Section#1{\s@ct@n{\kern-.3em #1}}

\def\Appendix{\setcounter{section}{0}\def\thesection{\Alph{section}}
 \let\th@s@ct@n\thesection
 \def\Section##1{\def\thesection{Appendix \Alph{section}}
 \s@ct@n{\kern-.3em##1}\let\thesection\th@s@ct@n
 \def\@currentlabel{\p@section\th@s@ct@n}}}

\newbox\commentbox

\newenvironment{abst}{\begingroup\narrower\noindent{\bf Abstract.}\enspace
 \ignorespaces}{\endgraf\endgroup}
\newenvironment{proof}{\begin{trivlist}\item{\it Proof.}\enspace\ignorespaces}
 {\qed\end{trivlist}}
\newenvironment{Proof}[1]{\begin{trivlist}\item{\it Proof of #1.}\enspace
 \ignorespaces}{\qed\end{trivlist}}
\newenvironment{example}[1]{\begin{trivlist}\item{\bf Example #1.}\enspace
 \ignorespaces}{\end{trivlist}}
\newenvironment{remark}{\begin{trivlist}\item{\sl Remark.}\enspace
 \ignorespaces}{\end{trivlist}}
 {\end{trivlist}\egroup\global\@ignoretrue\ifComments\unvbox\commentbox\fi
 \noindent}

\newtheorem{prop}{Proposition}[section]
\newtheorem{theor}[prop]{Theorem}
\newtheorem{lemma}[prop]{Lemma}
\newtheorem{corol}[prop]{Corollary}

\makeatother

\hybrid
\begin{document}

\begin{flushright}
\vp1
\vsk-3>
q-alg/9712002\hfill KYUSHU-MPS-1997-35\\
ITEP-TH-65/97\\
\vsk1.5>
\end{flushright}

\begin{center}
{\Large\bf On Solutions of the KZ and qKZ equations at Level Zero}\\
\vsk2>
{\large A. Nakayashiki$^{\,\star}$}
\vsk>
{\it Graduate School of Mathematics, Kyushu University,
Ropponmatsu 4-2-1, Fukuoka \,810, Japan}
\vsk1.6>
{\large S. Pakuliak$^{\,*}$}
\vsk>
{\it Bogoliubov Laboratory of Theoretical Physics, JINR,
Dubna, Moscow region \,141980, Russia}
\vsk.5>
{\sl and}
\vsk.5>
{\it Bogoliubov Institute of Theoretical Physics, Kiev \,252143, Ukraine}
\vsk1.6>
{\large V. Tarasov$^{\,\diamond}$}
\vsk>
{\it St.\;Petersburg Branch of Steklov Mathematical Institute,
Fontanka 27, St.\;Petersburg \,191011, Russia
\vsk.5>
{\sl and}
\vsk.5>
Department of Mathematics, Faculty of Science, Osaka University,
Toyonaka, Osaka \,560, Japan}
\end{center}
\ftext{\normalsize\sl$^\star$E-mail\/{\rm:} atsushi@rc.kyushu-u.ac.jp
\vsk0>\nt
$^*$E-mail\/{\rm:} pakuliak@thsun1.jinr.ru
\vsk0>\nt
$^\diamond${E-mail\/{\rm:} vt@pdmi.ras.ru}}
\vsk2>

\begin{abst}
We discuss relations between different formulae for solutions of
the Knizh\-nik-Za\-mo\-lod\-chi\-kov differential (KZ) and the quantum
Knizh\-nik-Za\-mo\-lod\-chi\-kov (qKZ) difference equations at level $0$
associated with rational solutions of the Yang-Baxter equation.
\end{abst}

\vsk>
\vsk0>

\setcounter{section}{0}
\setcounter{footnote}{0}

\Section{Introduction}
The quantum Knizhnik-Zamolodchikov (qKZ) difference equation has recently
attracted much attention. The qKZ equations appear naturally in
the representation theory of quantum affine algebras as equations for matrix
elements of products of vertex operators \cite{FR}. Later the qKZ equations
were derived as equations for traces of products of vertex operators
\cite{JM1}, their specializations being equations for correlation functions
in solvable lattice models. An important special case of the qKZ equations
was introduced earlier by Smirnov \cite{S1} as equations for form factors
in two-dimensional massive integrable models of quantum field theory.

In this paper we consider the rational qKZ equation associated with the Lie
algebra $\sltwo$. We will address the trigonometric qKZ equation in a separate
paper. We also restrict ourselves to the case of the qKZ equation at level
zero, which is precisely the case of Smirnov's equations for form factors.
Besides  this important application, this case is peculiar itself, which will
become clear in the paper. Another important special case is given by the qKZ
equation at level $-4$, which corresponds to equations for correlation
functions. We are going to discuss this case elsewhere. It is not much known
about solutions of the qKZ equation in the other cases than $\sltwo$. Only
a few results are available for the $\frak{sl}_N$ case \cite{S1}, \cite{KQ},
\cite{M}, \cite{TV1}.

There are several approaches to integral formulae for solutions of the qKZ
equation at level zero. The first one is given in \cite{S1}. Solutions
of the qKZ equation are obtained there in a rather starightforward way.
They are expressed via certain polynomials and they are enumerated by periodic
functions, which are arbitrary polynomials of exponentials of bounded degree.
Smirnov's construction can be seen as a deformation of hyperelliptic integrals,
the periodic functions being ``deformations'' of cycles on the corresponding
hyperelliptic curve \cite{S2}. From the view point of constructing solutions,
the main disadvantage of this approach is that it fails to work in the case of
the qKZ equation at nonzero level. As we understand now, the reason is that
Smirnov's formulae are intimately related to specific features of the qKZ
equation at level zero.

Another way to produce integral formulae for solutions of the qKZ equation is
given in \cite{JM1}, \cite{L}, \cite{KLP}. One has to calculate a trace of
a product of vertex operators over an infinite-dimensional representation of
the quantum affine algebra $U_q({\widehat{\sltwo})}$ or the centrally extended
Yangian double $\Aya$ \cite{L}, \cite{KLP}, using the bosonization technique,
thus getting integral solutions of the trigonometric and rational qKZ equation,
respectively. This approach is not restricted only to the qKZ equation at level
zero, but still cannot be applied to producing solutions of the qKZ equation in
general. The traces of products of vertex operators satisfy the qKZ equation
by construction and there is no need to verify independently that the final
integral formula gives a solution. On the other hand the resulting
parametrization of solutions by the products of vertex operators is very hard
to control effectively until now.

A general approach to integral representations for solutions of the qKZ
equation is developed in \cite{TV2}, \cite{TV3} combining ideas of \cite{SV},
\cite{V}, \cite{S1}. It allows to describe effectively the total space of
solutions of the qKZ equation at generic position, as well as to compute its
transition matrices between asymptotic solutions, which are substitutes for
monodromy matrices of a differential equation in the case of a difference
equation. Unfortunately, one cannot apply immediately these general results in
the case of the qKZ equation at level zero, because the genericity assumptions
imposed in \cite{TV1} are not fulfilled in this case. But it is possible to
repeat the consideration following the same lines as in \cite{TV1} making
necessary modifications and obtain a general formula for solutions of the qKZ
equation at level zero. As in Smirnov's formula, solutions are parametrized
by periodic functions which are polynomials in exponentials of bounded degree.
We beleive that Theorem \ref{main} describes all solutions of the qKZ equation
at level zero, though we have no proof of this conjecture yet.
For the recent development concerning general solutions of the qKZ equations
with values in finite-dimensional representantions and at rational levels
see \cite{MV}.

The general aim of this paper is to compare three above described types of
integral formulae. We will show that Smirnov's formula can be obtained from
a general formula \r{solutW} by a certain specialization of a periodic function
due to a certain trick available only at level zero. This trick was first
observed for the KZ differential equation. It turns out that the integrand in
a general integral formula for solutions of the KZ equation \cite{DJMM},
\cite{SV} becomes an exact form in the case in question, but it is possible
to produce nonzero solutions by integrating over suitable unclosed contours.
And a similar effect happens to occur in the difference case.

An open challenging problem is to describe traces of products of vertex
operators in terms of solutions of the qKZ equations given by the formula
\r{solutW}. Namely, for any product of vertex operators one must find
the corresponding periodic function. This will give a description of form
factors of local operators in an alternative way to Smirnov's axiomatic
approach \cite{S2}. In this paper we make this identification in the simplest
examples of the energy-momentum tensor and currents of the $SU(2)$-invariant
Thirring model.

The plan of the paper is as follows. In Section~2 we recall the definition of
the KZ and qKZ equations at level zero. We describe the necessary spaces of
rational functions in Section~3. Solutions of the KZ equation at level zero
are considered in Section~4. In Section~5 we describe the space of ``deformed
cycles''${\;=\;}$periodic functions, and the pairing between the spaces of
rational and periodic functions given by the hypergeometric integral.
We obtain solutions of the qKZ equation at level zero in Section~6.
The traces of products of vertex operators are considered in Section~7.

There are four Appendices in the paper which contain technical details and some
proofs.

\Section{The KZ and qKZ equations at level zero}
Let $\h$ be a nonzero complex number. Let ${V=\CC v_+ \oplus \CC v_-\,}$, and
$R(z)\in\hbox{End}\left(V^{\ot2}\right)$ be the following $R$-matrix:
\beq
\label{R-mat}
R(z)\,=\;{z+\h P\over z+\h},
\eeq
where $P\in\hbox{End}\left(V^{\ot2}\right)$ is the permutation operator.
We have
$$
R(z)\,=\,1+\h\>r(z)+O(\h^2)\,,\qqquad \h\to 0\,,
$$
where $r(z)$ is the corresponding classical $r$-matrix:
$$
r(z)\,=\,{-1+P\over z}\;.
$$

\goodbreak
Fix a nonzero complex number $p$ called {\em step\/}.
We consider the qKZ equation for a $V^{\ot n}$-valued function
$\psi(z_1\lc z_n)$:
\beq
\label{qKZ0}
\psi(z_1\lc z_j+p\lc z_n)\,=\,K_j(z_1\lc z_n)\>\psi(z_1\lc z_n)\,,
\eeq
$$
\vp{\Big|}
K_j(z_1\lc z_n)\,=\,R_{j,j-1}(z_j-z_{j-1}+p)\ldots R_{j,1}(z_j-z_1+p)
\>R_{j,n}(z_j-z_n)\ldots R_{j,j+1}(z_j-z_{j+1})\,,
$$
$j=1,\ldots n$. The number \,$-2+p/\h\,$ is called {\em level\/} of
the qKZ equation.
\Par
Let ${p=\h\)\ka}$ and consider the limit ${\h\to 0}$.
Then the qKZ difference equation turns into the KZ differential equation:
\beq
\label{KZ0}
\ka\,{\der\psi\over\der z_j}\;=\,
\sum_{\tsize{k=1\atop k\ne j}}^n r_{jk}(z_j-z_k)\,\psi\,,\qqquad j=1\lc n\,.
\eeq
The number \,$-2+\ka\,$ is called {\em level\/} of the KZ equation.
\vsk>
Let $\dsize\quad\sigma^+=\left({0\ \;1\atop 0\ \;0}\right)\,,\quad$
$\dsize\sigma^-=\left({0\ \;0\atop 1\ \;0}\right)\,,\quad$
$\dsize\sigma^3=\left({1\ \;\hp-0\atop 0\ \;{-1}}\right)\quad\;$ and
$\dsize\quad\Sigma^j=\tsize\Sum_{i=1}^n \sigma^j_i\,,\quad$ $j=\pm,3\,.\quad$
\vv>\\
The operators $\Sigma^\pm,\Sigma^3\in\hbox{End}\left(V^{\ot n}\right)$ defines
an $\sltwo$ action on $V^{\ot n}$.  For any $k$ such that ${0\le k\le n}$,
denote by $(V^{\ot n})_k$ the weight subspace
$$
(V^{\ot n})_k\,=\,\big\lbrace\,v\in V^{\ot n}\ \big|\ \Sigma^3 v=(n-2k)\,v\,
\big\rbrace.
$$
Notice that the operators $K_j(z_1\lc z_n)$, $j=1\lc n$, see \r{qKZ0},
commute with the $\frak{sl}_2$ action on $V^{\ot n}$.
\Par
In this paper we consider only the case of the KZ and qKZ equations
at level $0$, i.e.\ all over the paper we assume that
$$
\ka=2\qqquad\hbox{and}\qqquad p=2\h
$$
unless otherwise stated. Furthermore, fixing an integer $\ell$ such that
$0\leq 2\ell\leq n$, we discuss solutions of \r{qKZ0} and \r{KZ0} taking
values in $(V^{\ot n})_\ell$ and obeying an extra condition
\beq
\Sigma^+ \psi(z_1\lc z_n)\,=\,0\,,\vp{\Big|}
\label{sing}
\eeq
that is, the solutions taking values in singular vectors with respect to
the $\sltwo$ action.
\Par
In this paper we assume that $\Im\h<0$. Notice that we do not use this
assumption until the definition of the hypergeometic integral, see \r{pairing}.

\Section{The spaces of rational functions}
Let $M$ be a subset of $\{1\lc n\}$ such that \#$M=\ell$. We write
$M=\{m_1\lc m_\ell\}$ assuming that $m_1<\ldots<m_\ell$. The subset $M$
defines a point $\,\hat z_M\in\CC^\ell\,$:
\beq
\label{point}
\hat z_M=(z_{m_1}\lc z_{m_\ell})\,.
\eeq
For any subsets $M,N\subset\{1\lc n\}$ we say that $M\llee N$ if
$\# M=\# N$ and $m_i\leq n_i$ for any $i=1\lc\#M$.
For any function $f(t_1\lc t_\ell)$ we set
$$
\Asym f(t_1\lc t_\ell)\,=\sum_{\sigma\in\Sb_\ell}
\sgn(\sigma)\>f(t_{\sigma_1}\lc t_{\sigma_\ell})\,,
$$
and for any point $u=(u_1\lc u_\ell)\in \CC^\ell$ we define
$$
\Res{}\]f(u)\,=\,\res{}(\ldots\)\res{}
f(t_1\lc t_\ell)|_{t_\ell=u_\ell}\ldots)|_{t_1=u_1}\ .
$$

 From now on until the end of the section we assume that $z_1\lc z_n$ are
pairwise distinct complex numbers. Denote by $\Fc$ the space of rational
functions in $t$ with at most simple poles at points $z_1\lc z_n\,$. Let
$$
\Fc^{(k)}=\{f\in\Fc\ |\ f(t)=O(t^{k-2}),\ \ t\to\infty\},
\quad k\in\ZZ_{\geq0}\ .
$$
We consider $\Fol$ as a space of rational functions in $\ell$ variables
so that $\Fw\subset\Fol$ is the subspace of antisymmetric functions. Let
$$
\Fl = \{f\in\Fww\ |\res{\,t=z_m\!} f(t,t+\h,t_3\lc t_\ell)=0\,,
\quad m=1\lc n\}
$$
and
$$
\Fcl = \{f\in\Fww\ |\res{\,t_1=z_m\!}\,\Bigl(\,\res{\,t_2=z_m\!}\,
\Bigl(\dsize f(t_1,t_2,\lc t_\ell)\!\prod_{1\le a<b\le\ell}{1\over t_a-t_b}\,
\Bigr)\Bigr)=0\,,\quad m=1\lc n\}\,.
$$

For any $M=\{m_1,\ldots,m_\ell\}\subset\{1\lc n\}$ and $m\in M$,
\,let $\mu_M^{(m)}$, $g_M$, $w_M$ and $\tilde w_M$ be the following functions:
$$
\mu_M^{(m)}(t)\,=\,{1\over t-z_m}\prod _{l\in M\atop l\neq m}
\,{t-z_l-\h\over z_m-z_l-\h}\;,
$$
$$
\tilde w_M(t_1\lc t_\ell)\,=\,\Asym\Bigl(\,
\prod_{a=1}^\ell \mu_M^{(m_a)}(t_a)\Bigr)=\,
\det\bigl[\mu_M^{(m_a)}(t_b)\bigr]_{a,b=1}^\ell\,,
$$
$$
g_M(t_1\lc t_\ell)\,=\,\prod_{a=1}^\ell
\biggl({1\over t_a-z_{m_a}}\prod_{1\leq l<m_a}{t_a-z_l-\h\over t_a-z_l}\biggr)
\prod_{1\leq a<b\leq\ell}(t_a-t_b-\h)\,,
$$
\beq
\vp{\bigg|} w_M\,=\,\Asym\> g_M\ . \label{wM}
\eeq
The functions $\,\mu_M^{(m)cl}$, $g_M^{cl}$, $w_M^{cl}$ and $\tilde w_M^{cl}\,$
are defined by the same formulae at $\,\h=0\,$.

\begin{lemma}
\label{belong}
{\sl Let $M\subset\{1\lc n\}$, $\,\# M=\ell$. Then $\,w_M,\,\tilde w_M\in\Fl\,$
and $\,w_M^{cl},\,\tilde w_M^{cl}\in\Fcl\,$.}
\end{lemma}

\begin{lemma}
\label{res}
{\sl Let $M,N\subset\{1\lc n\}$, $\,\# M=\# N$. Then}
$\,{\Res{}w_M(\hat z_N)=0}\,$ {\sl unless $\,N\llee M$,}
$$
\Res{}\tilde w_M(\hat z_N)=\Res{}w_M^{cl}(\hat z_N)=
\Res{}\tilde w_M^{cl}(\hat z_N)=0\,,\qquad M\neq N\,,
$$
$$
\Res{}w_M(\hat z_M)=\Res{}g_M(\hat z_M)\,,\qqquad
\Res{}w_M^{cl}(\hat z_M)=\Res{}g_M^{cl}(\hat z_M)\,,
$$
$$
\Res{}\tilde w_M(\hat z_M)=\Res{}\tilde w_M^{cl}(\hat z_M)=1\,.
$$
\end{lemma}

\begin{prop}
\label{3}
{\sl Let $|z_j-z_k|\ne 0,\h\;$ for any $1\le j<k\le n\,$. Then the functions
$\,w_M$, $\,M\subset\{1\lc n\}$, $\,\# M=\ell$, \,form a basis in the space
$\Fl$, and the same do the functions $\,\tilde w_M$, $\,M\subset\{1\lc n\}$,
$\,\#M=\ell$.}
\end{prop}

\begin{prop}
\label{basis}
{\sl Let $\,z_1\lc z_n\,$ be pairwise distinct. Then the functions
$\,w_M^{cl}$, $\,M\subset\{1\lc n\}$, $\,\#M=\ell$, form a basis in the space
$\Fcl$, and the same do the functions $\,\tilde w_M^{cl}$,
$\,M\subset\{1\lc n\}$, $\,\# M=\ell$.}
\end{prop}
The propositions are proved in Appendix~\ref{AF}.
\Par
 From Lemma~\ref{res} and Propositions~\ref{3}, \ref{basis} it is clear that
\beq
w_M\,=\!\sum_{\scr N\subset\{1\lc n\}\atop \scr \# N=\# M}
\tilde w_N\,\Res{}w_M(\hat z_N)\,.
\label{w-tw}
\eeq
\beq
w_M^{cl}\,=\,\tilde w_M^{cl}\,\Res{}w_M^{cl}(\hat z_M)\,.
\label{w-twcl}
\eeq

\begin{lemma}\label{3*}
{\sl
For any $M\subset\{1\lc n\}$, $\# M=\ell-1$, the following relations hold:}
\bea
&\h\sum_{k\not\in M} w_{M\cup\{k\}}(t_1\lc t_\ell)\,={}
\label{wMk}
\\
&\quad {}=\,\Asym\!\!\left(\biggl(\;\prod_{a=2}^\ell\,(t_1-t_a-\h)\,-
\prod_{a=2}^\ell\,(t_1-t_a+\h)\,\prod_{m=1}^n{t_1-z_m-\h \over t_1-z_m}
\biggr)g_M(t_2\lc t_\ell)\right) .
\nn
\eea
\bea
&\sum_{k\not\in M} w_{M\cup\{k\}}^{cl}(t_1\lc t_\ell)\,={}
\label{wMkcl}
\\
&\quad {}=\,\Asym\!\!\left(\biggl(\;\sum_{m=1}^n\,{1\over t_1-z_m}\;-
\sum_{a=2}^\ell\,{2\over t_1-t_a}\biggr)\prod_{a=2}^\ell\,(t_1-t_a)\,
g_M^{cl}(t_2\lc t_\ell)\right) .
\nn
\eea
\end{lemma}
The Lemma is proved in Appendix~\ref{A3*}.
\Par
Let $M_{\rm ext}=\{1\lc\ell\}$. Say that $M_{\rm ext}$ is the {\em extremal}
subset. Due to Lemma~\ref{res} the right hand side of \r{w-tw} for
the extremal subset contains only one summand and
\beq
w_{M_{\rm ext}}(t_1\lc t_\ell)\,=
\prod_{1\leq a<b\leq \ell}\!{(z_a-z_b-\h)\>(z_a-z_b+\h)\over z_a-z_b}
\ \tilde w_{M_{\rm ext}}(t_1\lc t_\ell)\,.
\label{w-tw-ext}
\eeq

\Section{Solutions of the KZ equation at level zero}
Let $\;\dsize\phi_{cl}\,=\,\prod_{m=1}^n\,(t-z_m)^{-1/2}\,$ be
the {\em phase function\/}. Consider the hyperelliptic curve $\,\Ec\,$:
\beq
\smash{y^2\,=\,\prod_{m=1}^n\,(t-z_m)\,.}\vp{\prod_m}
\label{y2}
\eeq
Say that a contour $\gm$ on the curve $\Ec$ is {\em admissible\/} if
$\,\dsize\smash{\int_{\!\gm\,}\phi_{cl}\>f\>dt}\,$ is convergent and
$\,\dsize\smash{\int_{\!\gm\,}{\der\over\der t}(\phi_{cl}\>f)\>dt\,=\,0}\,$
for any $f\in\Fc^{(\ell+1)}\,$.

Denote by \,$\xi:\Ec\to\CC P^1\,$ the canonical projection:
$\,\xi:(t,y)\mapsto t\,$.
\begin{lemma}
\label{adm}
{\sl Let $\gm$ be a smooth simple contour on $\Ec$ avoiding the branching
points $z_1\lc z_n$. Assume that either $\gm$ is a cycle or $\gm$ admits
a parametrisation $\,\rho:\RR\to\Ec\)$, such that
$\,\bigl|\xi\bigl(\rho(u)\bigr)\bigr|\to\infty\,$ and
$\,\arg\xi\bigl(\rho(u)\bigr)\,$ has finite limits as $\,u\to\pm\infty\,$.
Then $\gm$ is admissible.}
\end{lemma}
\begin{proof}
If $\gm$ is a cycle, then the claim is clear. Since $\,2\ell\le n\,$ we have
that $\,\phi_{cl}(t)f(t)=O(t^{\ell-1-n/2})=O(t^{-1})\,$ as $\,t\to\infty\,$,
which proves the claim in the second case.
\end{proof}
We denote by $\,\Ch\,$ the set of all admissible contours and by
$\,\Cc\subset\Ch\,$ the set of cycles on $\Ec$.
With any subset $M\subset\{1\lc n\}$ we associate a vector
$v_M\in\left(V^{\ot n}\right)_\ell\,$ by the rule:
\beq
v_M=v_{\ep_1}\ot\cdots\ot v_{\ep_n}\,,
\label{vM}
\eeq
where $\,\ep_i=+\,$ for $\,i\not\in M\,$ and $\,\ep_i=-\,$ for $\,i\in M$.
For any $\,\gm_1\lc\gm_\ell\in\Ch\,$ we define a function
$\,\psi_{\gm_1\lc\gm_\ell}(z_1\lc z_n)\,$ as follows:
\beq
\psi_{\gm_1\lc\gm_\ell}(z_1\lc z_n)\,=\!\!
\sum_{\tsize{M\subset\{1\lc n\}\atop\# M=\ell}} v_M\!\!
\int_{\gm_1\times\ldots\times\gm_\ell}\!\!\!w_M^{cl}(t_1\lc t_\ell)\,
\prod_{a=1}^{\ell}\phi_{cl}(t_a)\;dt_1\wedge\ldots\wedge dt_\ell\,.
\label{psiSV}
\eeq
\begin{theor}\label{SV}
{\rm \cite{DJMM}, \cite{SV}}\enspace
{\sl For any $\,\gm_1\lc\gm_\ell\in\Ch\,$ the function
$\,\psi_{\gm_1\lc\gm_\ell}(z_1\lc z_n)\,$
\vv.2>
is a solution of the KZ equation\/} {\rm\r{KZ0}} {\sl taking values in
$(V^{\ot n})_\ell$ and satisfying the condition\/} {\rm\r{sing}}.
\end{theor}
\begin{remark}
Formula \r{psiSV} gives a solution of the KZ equation for arbitrary $n,\,\ell$.
For $2\ell=n$, there exists another integral representation for solutions of
the KZ equation, see the formula \r{psiS} and Theorem~\ref{Sm}, which is proved
in \cite{N1} for a general $\frak{sl}_N$ case in a straightforward way.
This integral representation is a quasiclassical limit of Smirnov's formulae
for solutions of the qKZ equation. We will show that the formula \r{psiS} is
a specialization of the formula \r{psiSV} for a certain choice of
the integration contours, cf.\ \r{psiSVS}.
\end{remark}

Assume until the end of this section that $2\ell=n$. This means that we are
looking for a singlet solution $\psi(z_1\lc z_n)$ of the KZ equation \r{KZ0}:
$$
\Sigma^j\psi(z_1\lc z_n)=0\,,\qqquad j=\pm,3\,.
$$
For any $M\subset\{1\lc n\}$, \,$\# M=\ell$, define functions
$\nu_M$ and $\tilde\nu_M$ by the formulae:
$$
\nu_M^{cl}(t)= \prod_{j\not\in M}(t-z_j),\qqquad
\tilde\nu_M^{cl}(t)=\nu_M^{cl}(t)\,\Bigl(\,\sum_{j\in M}{1\over{t-z_j}}\;-
\sum_{j\not\in M}{1\over{t-z_j}}\,\Bigr)\,.
$$
It is clear that
\beq
\phi_{cl}\,\tilde\nu_M^{cl}\,=\,
-2\,{d\over dt}\bigl(\phi_{cl}\,\nu_M^{cl}\bigr)\,.
\label{totcl}
\eeq

Since $n$ is even, the curve $\Ec$ has no branching point at infinity.
We distinguish two infinity points $\infty^\pm\in\Ec$ as follows:
$$
y=\pm\,t^{n/2}\,\bigl(1+o(1)\bigr)\,,\qqquad (t,y)\to\infty^\pm\,.
$$
Fix an admissible contour $\gm^\infty$ going from $\infty^+$ to $\infty^-$,
cf.\ Lemma~\ref{adm}.
\begin{lemma}
\label{gm8}
{\sl Let $2\ell=n$. Then}
$\,\dsize\int_{\gm^\infty}\!\phi_{cl}\,\tilde\nu_M^{cl}\,dt\,=\,4\,$.
\end{lemma}
\begin{proof}
The statement immediately follows from the formula \r{totcl} and the definition
of $\gm^\infty$.
\end{proof}
\begin{lemma}
\label{7**cl}
{\sl Let $\,2\ell=n$. Then, for any $\,M\subset\{1\lc n\}$, $\,\#M=\ell$,
we have}
$$
\tilde\nu_M^{cl}(t)=\sum_{m\in M} \mu_M^{(m)cl}(t)\,
\res{}\tilde\nu_M^{cl}(z_m)\,.
$$
\end{lemma}
\begin{proof}
The claim follows from Lemma~\ref{7**} since, using explicit formulae,
it is easy to see that
$$
\h^{-1}\tilde\nu_M\to\,\tilde\nu_M^{cl}\,,\qqquad
\h^{-1}\res{}\tilde\nu_M(z_m)\to\res{}\tilde\nu_M^{cl}(z_m)\,,\qqquad
\mu_M^{(m)}\to\,\mu_M^{(m)cl}\,,
$$
as $\,\h\to 0\,$.
\end{proof}
\begin{corol}
\label{twnu}
{\sl Let $2\ell=n$. Then, for any $M\subset\{1\lc n\}$, $\# M=\ell$, we have}
$$
\tilde w_M^{cl}(t_1\lc t_\ell)\,=\;{1\over\res{}\tilde\nu_M^{cl}(z_{m_\ell})}
\;\Asym\bigl(\mu_M^{(m_1)cl}(t_1)\ldots\mu_M^{(m_{\ell-1})cl}(t_{\ell-1})
\tilde\nu_M^{cl}(t_\ell)\bigr)\,.
$$
\end{corol}

For any $\,\gm_1\lc\gm_{\ell-1}\in\Cc\,$ we define a function
$\,\psi^{S}_{\gm_1\lc\gm_{\ell-1}}(z_1\lc z_n)\,$:
\bea
&\kern-20pt& \psi^{S}_{\gm_1\lc\gm_{\ell-1}}(z_1\lc z_n)\,={}
\label{psiS}
\\
&\kern-20pt&\quad{}=\!\!\sum_{\tsize{M\subset\{1\lc n\}\atop\# M=\ell}}
\prod_{1\le a<b\le\ell}\!(z_{m_a}-z_{m_b})\,
\prod_{k\not\in M}(z_{m_\ell}-z_k)^{-1}\,\det\Bigl[\,
\int_{\gm_b}\phi^{cl}\>\mu_M^{(m_a)cl}\>dt\,\Bigr]_{a,b=1}^{\ell-1}\;v_M\,.
\nn
\eea

\begin{theor}\label{Sm}
{\rm \cite{N1}}\enspace
{\sl For any $\,\gm_1\lc\gm_{\ell-1}\in\Cc\,$ the function
$\,\psi^S_{\gm_1\lc\gm_{\ell-1}}(z_1\lc z_n)\,$
\vv.2>
is a solution of the KZ equation\/} {\rm\r{KZ0}} {\sl taking values in
$(V^{\ot n})_\ell$ and satisfying condition\/} {\rm\r{sing}}.
\end{theor}
\begin{proof}
Subsequently using formulae \r{psiSV}, \r{w-twcl}, Corollary~\ref{twnu},
Lemma~\ref{gm8} and Lemma~\ref{res} we obtain
\bea
&\kern-4em&
\psi_{\gm_1\lc\gm_{\ell-1},\gm^\infty}(z_1\lc z_n) \,=
\sum_M\; v_M\,\Res{}w_M^{cl}(\hat z_M)\;\det\Bigl[\,
\int_{\gm_b}\phi^{cl}\>\mu_M^{(m_a)cl}\>dt\,\Bigr]_{a,b=1}^\ell
\label{psiSVS}
\\
&\kern-4em&\quad {}\,=\,4\,\sum_M\; v_M\;
{\Res{}w_M^{cl}(\hat z_M)\over\res{}\tilde\nu_M^{cl}(z_{m_\ell})}
\ \det\Bigl[\,
\int_{\gm_b}\phi^{cl}\>\mu_M^{(m_a)cl}\>dt\,\Bigr]_{a,b=1}^{\ell-1}
\,=\,4\;\psi^{S}_{\gm_1\lc\gm_{\ell-1}}(z_1\lc z_n)\,.
\nn
\eea
Here $\gm_\ell=\gm^\infty$. Since $\,\gm_1\lc\gm_{\ell-1},\gm^\infty\in\Ch\,$,
the theorem follows from Theorem~\ref{SV}.
\end{proof}
\begin{remark}
Formula \r{totcl} and Corollary \ref{twnu} show that for $n=2\ell$, that is
in the zero weight case, the integrand in the right hand side of the formula
\r{psiSV} for solutions of the KZ equation is an exact form. So that, for
any $\,\gm_1\lc\gm_\ell\in\Cc\,$ we have
$$
\psi^S_{\gm_1\lc\gm_\ell}(z_1\lc z_n)\,=\,0\,.
$$
Therefore, to produce nonzero solutions of the qKZ equation we inevitably have
to consider not only cycles on the curve $\Ec$, but also certain unclosed
contours, which we call admissible.
\end{remark}

\Section{The space of ``deformed cycles'' and the hypergeometric integral}
Let $\Fq$ be the space of functions $F(t;z_1\lc z_n)$ such that
\beq
F(t;z_1\lc z_n)\,\prod_{j=1}^n\,(1-\ee^{2\pi i(t-z_j)/p})
\,=\,P\bigl(\ee^{2\pi it/p}\bigr)
\label{cycles}
\eeq
is a polynomial in $\ee^{2\pi i t/p}$ of degree at most $\,n\,$ and
$$
F(t;z_1\lc z_j+p\lc z_n)\,=\,F(t;z_1\lc z_n),\qquad j=1\lc n\,.
$$
Notice that the definition \r{cycles} implies that any $F(t)\in\Fq$ is
a periodic function of $t$:\ \ $F(t+p)=F(t)$.

For any function $F\in\Fq$ we set
\beq
F(\pm\infty\);z_1\lc z_n)\,=\lim{t\to\pm\infty}\]F(t;z_1\lc z_n)\,.
\eeq

Let $\Fcq\subset\Fq$ be the subspace of functions $F(t;z_1\lc z_n)$ such
that $F(\pm\infty\);z_1\lc z_n)=0$, that is the polynomial $P$ in \r{cycles}
obeys extra conditions
$$
P(0)=0\qquad\hbox{\ and}\qquad\deg P<n\,.
$$

Let $\Qc\subset\Fq$ be the following subspace: $\Qc=\CC\,1\oplus\CC\,\Theta$,
where
\beq
\Theta(t)\,=\,\prod_{j=1}^n\,
{1+\ee^{2\pi i(t-z_j)/p}\over 1-\ee^{2\pi i(t-z_j)/p}}\;.
\label{Theta}
\eeq

Let $\phi(t)$ be the {\em phase function\/}:
\beq
\phi(t)\,=\,\prod_{j=1}^n\,
{\Ga{\fract{t-z_j-\h}{p}}\over \Ga{\fract{t-z_j}{p}} }\;.
\label{phase}
\eeq
\begin{lemma}
\label{asympt}
{\sl For any $\,\ep>0\,$ the phase function $\phi(t)$ has the following
asymptotics
$$
\phi(t)\,=\,(t/p)^{-n\h/p}\bigl(1+o(1)\bigr)
$$
as \,$t\to\infty$, \ $\varepsilon<|\arg\,(t/p)|<\pi-\varepsilon\,$.}
\end{lemma}
\nt
The statement follows from the Stirling formula.
\Par
Denote by $D$ the operator defined by
\beq
D\]f(t)\,=\,f(t)-f(t+p)\>\prod_{j=1}^n {t-z_j-\h\over t-z_j}\;,
\label{tot-dif}
\eeq
which means
$$
\phi(t)\,D\]f(t)\,=\,\phi(t)\>f(t)-\phi(t+p)\>f(t+p)\,.
$$
We call the functions of the form $D\]f$ the {\em total differences\/}.
For any $s\in\ZZ\,$ let $\,\Ucup_s^\pm$ be the following sets:
\beq
\Ucup^+_s=\Cup_{j=1}^n(z_j+\h+p\)\ZZ_{\le s})\,,\qquad\qquad
\Ucup^-_s=\Cup_{j=1}^n(z_j-p\)\ZZ_{\le s})\,.
\label{UU}
\eeq
Let $I(w,W)$ be the {\em hypergeometric integral\/}:
\beq
I(w,W)=\int_{C}\phi\>w\)W\>dt\,,
\label{pairing}
\eeq
where $C$ is a simple curve separating the sets $\,\Ucup^+_1$ and
$\,\Ucup^-_1$, and going from $-\infty$ to $+\infty$. More precisely,
the contour $C$ admits a parametrization $\rho:\RR\to\CC$ such that
$\rho(u)\to\pm\infty$ and $\Im\rho(u)$ has finite limits as $u\to\pm\infty$.
Recall that we assume $\,\Im p<0\,$ and $\,p=2\h\,$.
\begin{remark}
A construction of the hypergeometric integral for the qKZ equation at general
level is given in \cite{TV2} for the rational case and in \cite{TV3} for
the trigonometric case. In this paper we adapt the general construction from
\cite{TV2} to the case of the qKZ equation at level zero.
\end{remark}
\begin{lemma}\label{4*}
{\sl Let ${w\in\Fc,\ W\in\Fcq}\;$ or $\,{w\in\Fll,\ W\in\Fq}$.
Then the integral $I(w,W)$ is ab\-so\-lu\-te\-ly
con\-ver\-gent and does not depend on
a particular choice of the contour $C$.}
\end{lemma}

\begin{proof}
The integrand $\phi(t)\>w(t)\)W(t)$ of the integral $I(w,W)$ behaves like
$O(t^{-2})$ as $t\to\pm\infty$, which proves the convergence of the integral.

The poles of the integrand ${\phi\>w\)W}$ belong to the set
$\Ucup^+_0\!\cup\)\Ucup^-_0$. Therefore, a homotopy class of the contour $C$
in the complement of the singularities of the integrand does not depend on
a particular choice of the contour and the same does the value of the integral
$I(w,W)$.
\end{proof}

\begin{lemma}\label{44}
{\sl Assume that $\;{w\in\Fll,\ W\in\Qc}$. Then we have $I(w,W)=0$.}
\end{lemma}

\begin{proof}
Let $W=1\,$. Then the integrand of the integral $I(w,1)$ equals $\phi(t)\>w(t)$
and has no poles at points of the set $\,\Ucup^-_1$. Moreover, the integrand
uniformly behaves like $O(t^{-2})$ as $t\to\infty$ in the semiplane
$\Im t\le 0$. Therefore, for any large negative $A$, the contour $C$ can be
replaced by the line $\{t\in\CC\ |\ \Im t=A\}$ without changing the integral
$I(w,1)$. Tending $A$ to $-\infty$ we obtain that $I(w,1)=0$.

If $W=\Theta\,$, then the integrand $\phi(t)\>w(t)\>\Theta(t)$ has no
poles at points of the set $\,\Ucup^+_1$ and uniformly behaves like $O(t^{-2})$
as $t\to\infty$ in the semiplane $\Im t\ge 0$. Therefore, for any large
positive $A$, the contour $C$ can be replaced by the line
$\{t\in\CC\ |\ \Im t=A\}$ without changing the integral $I(w,\Theta)$,
and tending $A$ to $+\infty$ we obtain that $I(w,\Theta)=0$.
\end{proof}

\begin{lemma}\label{4}
{\sl Assume that ${w\in D\Fc,\ W\in\Fcq}\;$ or
$\;{w\in D\Fc^{(\ell+1)},\ W\in\Fq}\;$. Then we have $I(w,W)=0$.}
\end{lemma}

\begin{proof}
Let $C$ be the contour in the formula \r{pairing}. For any $A>0$ set
${C^A=\,\{\,t\in C\ |\ |\Re t|\le A\,\}}$\,. Let $C^A_p$ be the contour $C^A$
shifted by $p\,$: \,$t\in C^A\,\Longleftrightarrow\,(t+p)\in C^A_p$. Let
$\,\Dl^A_\pm\,$ be segments of length $|p|$ such that the contour
$\,C^A+\Dl^A_+-C^A_p-\Dl^A_-\,$ is closed and $\,\pm\Re t>0\,$ for
$t\in\Dl^A_\pm$.

Let $w=D\tilde w$. The integral defining $I(w,W)$ is convergent. Thus we have
$$
I(w,W)\,=\>\lim{A\to\infty}\>\int_{C^A}\phi\>w\)W\>dt\,.
$$
Using formulae \r{tot-dif} we obtain
$$
\int_{C^A}\phi\>w\)W\>dt\ =\int_{C^A-C^A_p}\!\!\!\phi\>\tilde w\)W\>dt\ =
\int_{C^A+\Dl^A_+-\>C^A_p-\Dl^A_-}\!\!\!\!\!\!\!\!\phi\>\tilde w\)W\>dt
\ +\int_{\Dl^A_- -\Dl^A_+}\!\!\!\phi\>\tilde w\)W\>dt\,.
$$
Since there is no poles of the integrand $\,{\phi\>\tilde w\)W}\,$ inside
the contour $\,{C^A-\Dl^A_+-C^A_p+\Dl^A_-}\,$, \,the cor\-res\-pond\-ing
integral equals zero. Moreover,
$$
\lim{A\to\infty}\!\int_{\Dl^A_- -\Dl^A_+}\!\!\!\phi\>\tilde w\)W\>dt\,=
\>p\]\lim{A\to\infty}\bigl(\phi(-A)\>\tilde w(-A)\)W(-A)-
\phi(A)\>\tilde w(A)\)W(A)\bigr)\,=\,0
$$
under the assumptions of the lemma. The lemma is proved.
\end{proof}

\begin{lemma}
\label{48}
{\sl Let a function $f$ be such that $\,{(f-D\]f)\in\Fc^{(\ell+1)}}$.
Then for any $\,{W\in\Fq}\,$ we have that $\,I(D\]f,W)=0$.}
\end{lemma}
The proof is similar to the proof of the second part of Lemma~\ref{4}.

\begin{remark}
Lemmas~\ref{4*}, \ref{44} and \ref{7*} below remain valid under weaker
assumptions that the contour $C$ separates the sets $\,\Ucup^+_0$ and
$\,\Ucup^-_0$. For Lemma~\ref{4} it suffices to assume that $C$ separates
the sets $\,\Ucup^+_0$ and $\,\Ucup^-_1$, and for Lemma~\ref{48}, that $C$
separates the sets $\,\Ucup^+_1$ and $\,\Ucup^-_0$.
\end{remark}

\Section{Solutions of the qKZ equation at level $0$}
\label{666}
Denote by $I^{\ot \ell}(w,W)$ the following integral:
\beq
I^{\ot\ell}(w,W)=\int_{C^\ell}w(t_1\lc t_\ell)\>W(t_1\lc t_\ell)\,
\prod_{a=1}^\ell \phi(t_a)\,dt_a\,.
\label{pair-ell}
\eeq

For any $W\in \Fq^{\ot\ell}$, define a function $\psi_W(z_1\lc z_n)$ with
values in $(V^{\ot n})_\ell$ as follows:
\beq
\psi_W(z_1\lc z_n)=
\sum_{\tsize{M\subset\{1\lc n\}\atop\# M=\ell}}\!\!\!I^{\ot\ell}(w_M,W)\;v_M\,,
\label{solutW}
\eeq
where $\,v_M\,$ is defined by \r{vM}.
\begin{remark}
It is clear that for any $\,W\in\Fq^{\ot\ell}\,$ we have
$\dsize\,\psi_W={1\over\ell\>!}\,\psi_{\Asym\!W}$.
\end{remark}
\begin{prop}\label{5*}
{\sl For any \,$W\in\Qc\ot\Fq^{\ot(\ell-1)}$ we have
\,\,$\psi_W(z_1\lc z_n)=0$\,.}
\end{prop}
The claim follows from the Lemma~\ref{44}.

\begin{prop}\label{5}
{\sl For any $W\in\Fq^{\ot\ell}$ the function $\psi_W(z_1\lc z_n)$
satisfies the condition}
$$
\Sigma^+\psi_W(z_1\lc z_n)=0\,.
$$
\end{prop}
The claim follows from the Lemmas~\ref{3*} and \ref{4}.

\begin{theor}
\label{main}
{\sl For any $W\in\Fq^{\ot\ell}$ the function $\psi_W(z_1\lc z_n)$
is a solution of the qKZ equation} {\rm\r{qKZ0}} {\sl taking values in
$(V^{\ot n})_{\ell}$}.
\end{theor}
\begin{proof}
The statement follows from the results on formal integral representations
for solutions of the qKZ equation \cite{V} and Lemmas~\ref{4}, \ref{48}.
We give more details in Appendix~\ref{proofmain}.
\end{proof}

Let
\beq
\tilde v_M=\sum_{\scr N\subset\{1\lc n\}\atop\scr \#N=\ell}\!\!\!v_N
\,\Res{}w_N(\hat z_M)\,,\qqquad M\subset\{1\lc n\}\,,\quad \#M=\ell\,,
\label{S-basis}
\eeq
be another basis in $(V^{\ot n})_\ell$. Then we have
\beq
\psi_W(z_1\lc z_n)=\sum_{\scr M\subset\{1\lc n\}\atop\scr\# M=\ell}
\!\!\!I^{\ot\ell}(\tilde w_M,W)\;\tilde v_M\ .
\label{solWSbas}
\eeq
For ${W=\prod\limits_{a=1}^n W_a}\,$, $\;{W_a\in\Fq}$, $\;{a=1\lc\ell}$,
\,the last formula can be written in the determinant form
\beq
\psi_W(z_1\lc z_n)=\sum_{\scr M\subset\{1\lc n\}\atop\scr\# M=\ell}\!\!\!
\det\bigl[ I(\mu_M^{(m_a)},W_b)\bigr]_{a,b=1}^\ell\;\tilde v_M\ .
\label{solWdet}
\eeq
The solutions $\psi_W(z_1\lc z_n)$ can be written also via suitable
polynomials rather than rational functions.
For any $\,{M\subset\{1\lc n\}}\,$ set
\beq
\label{PMP}
P^+_M(t)\,=\,\prod_{m\in M}\,(t-z_m-2\h)\,,\qqquad
P^-_M(t)\,=\,\prod_{k\not\in M}\,(t-z_k-2\h)\,.
\eeq
Denote by $\,T_\h$ the operator defined by $\;T_\h\]f(t)\,=\,f(t)-f(t+\h)\,$.
For any rational function $f(t)\>$ let $\,[f(t)]_+\in\CC[t]\,$ be its
polynomial part
$$
f(t)\,=\,[f(t)]_+\,+\,o(1)\,,\qqquad t\to\infty\,.
$$
Let $\,Q_M^{(1)}\lc Q_M^{(\ell)}\,$ be the polynomials given by
\beq
\label{QM}
Q_M^{(a)}(t)\,=\,P^+_M(t+h)\>\Bigl[\>T_\h\Bigl({P^-_M(t)\over P^+_M(t+h)}\,
\Bigl[{P^+_M(t+h)\over(t+\h)^a}\Bigr]_+\Bigr)\Bigr]_+\;+\;
P^-_M(t)\>\Bigl[\>T_\h\Bigl({P^+_M(t)\over t^a}\Bigr)\Bigr]_+\,.
\eeq
\begin{remark}
The polynomials $\,Q_M^{(a)}\,$ coincide with the polynomials introduced
in \cite{S1} modulo a certain change of variables. Let us write the dependence
of the polynomials $\,Q_M^{(a)}\,$ on $z_1\lc z_n$ explicitly:
$Q_M^{(a)}(t)\)=\)Q_M^{(a)}(t;z_1\lc z_n)\,$.
Let $A^c_a(t\vert \la_1\lc\la_{\ell}\vert \mu_1\lc\mu_{\ell-2c})$
be the polynomial defined by (118) in \cite{S1}. Let $c=0$ or $c=-1$ and
$n=2\ell-2c$. Then
$$
A^c_a(t\vert \hat z_M \vert \hat z_\Mbar)\,=\,
Q^{(\ell-c-a)}_M(t-\h/2;z_1-2\h\lc z_n-2\h)\,,\qquad a=1\lc\ell-c-1,
$$ where $\Mbar=\{1\lc n\}\setminus\]M$ and $\hat{z}_{\Mbar}$ is defined
similarly to $\hat z_M$, cf.\ \r{point}. Notice that there is a misprint in
(118) in \cite{S1}.

The polynomials $\,Q_M^{(a)}\,$ are the rational analogues of the polynomials
$A^{(j,k)}_b$ introduced in \cite{JKMQ} for the trigonometric case.
The precise correspondence is as follows:
$$
A^{(\ell,n-\ell)}_a(t\vert\hat{z}_M\vert\hat{z}_\Mbar)\,\sim\,
Q^{(\ell+1-a)}_M(t-\h/2;z_1-2\h\lc z_n-2\h)\,,\qqquad a=1\lc\ell\,.
$$
\end{remark}
\vsk>
\begin{lemma}\label{7procent}
{\sl For any ${M\subset\{1\lc n\}}$, \,${\#M=\ell}$, and any
${m\in M}$ the following identity holds\/}{\rm:}
$$
\kern-10pt
D\>\Bigl(\prod_{\tsize{k=1\atop k\neq m}}^n\,(t-z_k-2\h)\Bigr)\,=\,
\h\prod_{\tsize{k=1\atop k\neq m}}^n\,(z_m-z_k-\h)\>\mu_M^{(m)}(t)\,+\,
\sum_{a=1}^{\ell}\,Q_M^{(a)}(t)\>(z_m+2\h)^{a-1}\,.
$$
\end{lemma}
The lemma is proved in Appendix~\ref{DQ}.
\goodbreak
\Par
For any $\,{M\subset\{1\lc n\}}\,$ set
\beq
\label{S-basis2}
v^S_M\,=\prod_{1\leq a<b\leq \ell}\,
{(z_{m_a}\!-z_{m_b})\over(z_{m_a}\!-z_{m_b}-\h)\>(z_{m_a}\!-z_{m_b}+\h)}
\ \tilde v_M\,,
\eeq
cf.\ \r{w-tw-ext}, \r{S-basis}, which provides that
$\;{v^S_{M_{\rm ext}}\!=\,v_{M_{\rm ext}}^{\vp1}\!+\ldots}\;$ where dots
stand for a linear combination of vectors $\,v_M$, $\,M\ne M_{\rm ext}$,
$\,\#M=\ell$. The main property of this basis is given by
Corollary \ref{symbasis}.
\Par
Using Lemmas~\ref{4},\,\ref{7procent} and formulae \r{solWdet}, \r{S-basis2}
we can rewrite the solution $\psi_W(z_1\lc z_n)$, $W\in\Fcq^{\ot\ell}$,
via the polynomials $\,Q_M^{(1)}\lc Q_M^{(\ell)}\,$ and the basis
$\{v^S_M\}_{\#M=\ell}\,$. Let $\,W_1\lc W_\ell\in\Fcq$ \,and
$\,W(t_1\lc t_\ell)=W_1(t_1)\ldots W_\ell(t_\ell)$. Then
\beq
\psi_W(z_1\lc z_n)=(-\h)^{-\ell}
\sum_{\tsize{M\subset\{1\lc n\}\atop\# M=\ell}}
\prod_{\tsize{k\not\in M\atop m\in M}}{1\over z_m-z_k-\h}\;
\det\bigl[I(Q_M^{(a)},W_b)\bigr]_{a,b=1}^\ell\;v^S_M\,.
\label{det-ar-spin}
\eeq
Observe that if $Q$ is a polynomial and ${W\in\Fcq}$, then the integrand of
the integral $I(Q,W)$ has no poles at points $z_m,\,z_m\!-p,\,z_m\!+\h+p$,
$\;m=1\lc n$, and we have
\beq
I(Q,W)=\int_{C'}\phi\,Q\>W\>dt
\label{IQW}
\eeq
for any contour $C'$ going from $-\infty$ to $+\infty\,$ and separating
the sets $\,\Ucup^+_0$ and $\,\Ucup^-_{-1}$. For instance,
if $\,z_1\lc z_n\,$ are real, then we can take $C'=3\h/2+\RR$.
\Par
Let $2\ell=n$. This means that we consider a singlet solution
$\psi(z_1\lc z_n)$ of the qKZ equation \r{qKZ0}:
$$
\Sigma^j\psi(z_1\lc z_n)=0\,,\qqquad j=\pm,3\,.
$$
For any $M\subset\{1\lc n\}$, \,$\# M=\ell$, \,let $\,\nu_M\,$ and
$\,\tilde\nu_M\,$ be the following functions:
$$
\nu_M(t)= \prod_{k\not\in M}(t-z_k-2\h),\qquad
\tilde\nu_M(t)=\prod_{k\not\in M}(t-z_k-2\h)
-\prod_{k\not\in M}(t-z_k-\h)\prod_{m\in M} {(t-z_m-\h)\over(t-z_m)}\ .
$$
It is clear that
\beq
\tilde\nu_M=D\nu_M\ .
\label{total}
\eeq
\begin{lemma}
\label{7*}
{\sl Let $2\ell=n$. Then
$\;\dis I(\tilde\nu_M,W)=\,p^{\ell+1}\>
\bigl(W(-\infty)-W(+\infty)\bigr)\;$ for any $W\in\Fq$.}
\end{lemma}

\begin{proof}
We use notations from the proof of Lemma~\ref{4}.
The integrand $\phi(t)\>\tilde\nu_M(t)\)W(t)$ behaves like $O(t^{-2})$ as
${t\to\pm\infty}$. Hence, the integral $I(\tilde\nu_M,W)$ is convergent
and we have
$$
I(\tilde\nu_M,W)\,=\lim{A\to\infty}\>\int_{C^A}\phi\>\tilde\nu_M W\>dt\,.
$$
Similarly to the proof of Lemma~\ref{4}, using the formula \r{total}, we obtain
\bea
\kern-2em
I(\tilde\nu_M,W)\,={}&\kern-18pt&
\lim{A\to\infty}\!\int_{\Dl^A_- -\Dl^A_+}\!\!\!\phi\>\nu_M W\>dt\;={}
\nn
\\
{}={}&\kern-18pt& p\]\lim{A\to\infty}\bigl(\phi(-A)\>\nu_M(-A)\)W(-A)-
\phi(A)\>\nu_M(A)\)W(A)\bigr)\,=\,
p^{\ell+1}\>\bigl(W(-\infty)-W(+\infty)\bigr)\,.
\nn
\eea
The lemma is proved.
\end{proof}

\begin{lemma}\label{7**}
{\sl Let $2\ell=n$. For any $M\subset\{1\lc n\}$, $\# M=\ell$, we have}
$$
\tilde\nu_M(t)\,=\sum_{m\in M}\>\mu_M^{(m)}(t)\;\res{}\tilde\nu_M(z_m)\ .
$$
\end{lemma}
\begin{proof}
Both sides of the formula are rational functions in $t$ with at most simple
poles at points $z_m$, $\,m\in M$, and have the same growth $O(t^{\ell-2})$ as
$\,t\to\infty$. Moreover, they have the same residues at points $z_m$,
$\,m\in M$, and the same values at the points $z_m+\h$, $\,m\in M$, which
completes the proof.
\end{proof}
\begin{corol}\label{77}
{\sl Let $2\ell=n$. Then for any $M\subset\{1\lc n\}$, $\# M=\ell$, we have}
$$
\tilde w_M(t_1\lc t_\ell)\,=\;
{1\over\res{}\tilde\nu_M(z_{m_\ell})}\;\Asym\bigl(\mu_M^{(m_1)}(t_1)\ldots
\mu_M^{(m_{\ell-1})}(t_{\ell-1})\tilde\nu_M(t_\ell)\bigr)\ .
$$
\end{corol}
The last corollary, the formula \r{solWSbas} and Lemma~\ref{7*} imply the next
theorem.

\begin{theor}
\label{7@}
{\sl Let $2\ell=n$. Let $\,W_1\lc W_{\ell-1}\in\Fcq$, $\,W_\ell\in\Fq$ and
$\,W(t_1\lc t_\ell)=W_1(t_1)\ldots W_\ell(t_\ell)$. Then}
$$
\psi_W(z_1\lc z_n)\,=\,p^{\ell+1}\>
\bigl(W_\ell(+\infty)-W_\ell(-\infty)\bigr)\!\!
\sum_{\tsize{M\subset\{1\lc n\}\atop\scr \# M=\ell}}
\,{\prod\limits_{a=1}^{\ell-1}\,(z_{m_\ell}-z_{m_a})\over
\!\prod\limits_{k=1}^n\>(z_{m_\ell}-z_k-\h)}\;
\det\bigl[I(\mu_M^{(m_a)},W_b\bigr]_{a,b=1}^{\ell-1}\ \tilde v_M\,.
$$
\end{theor}

Now we can rewrite the solution of the qKZ equation at level $0$ given by
Theorem~\ref{7@} via the polynomials $\,Q_M^{(1)}\lc Q_M^{(\ell)}\,$ and
recover the formulae from \cite{S1} for the singlet form-factors in
the $SU(2)$-invariant Thirring model. For $2\ell=n$ the polynomials can be
written in the form which appeared for the first time in the book \cite{S1}:
\beq
Q_M^{(a)}(t)\,=\,
P^+_M(t+h)\>\bigl[\>T_\h\bigl((t+\h)^{-a}\>P^-_M(t)\bigr)\bigr]_+\,+\,
P^-_M(t)\>\bigl[\>T_\h\bigl(t^{-a}\>P^+_M(t)\bigr)\bigr]_+\,,
\eeq
cf.\ \r{PMP}, \r{QM}.
Notice that the polynomial $\,Q_M^{(\ell)}\)$ vanishes identically
in this case.
\Par
Finally, using Lemmas~\ref{4},\,\ref{7procent} and formulae \r{solWSbas},
\r{S-basis2}, we have the following statement.

\begin{theor}
\label{7final}
{\sl Let $2\ell=n$. Let $\,W_1\lc W_{\ell-1}\in\Fcq$, $\,W_\ell\in\Fq$ and
$\,W(t_1\lc t_\ell)=W_1(t_1)\ldots W_\ell(t_\ell)$. Then}
$$
\psi_W(z_1\lc z_n)\,=\,2^\ell p\>
\bigl(W_\ell(+\infty)-W_\ell(-\infty)\bigr)
\!\!\sum_{\tsize{M\subset\{1\lc n\}\atop\#M=\ell}}
\prod_{\tsize{k\not\in M\atop m\in M}}{1\over z_k-z_m+\h}\;
\det\bigl[I(Q_M^{(a)},W_b)\bigr]_{a,b=1}^{\ell-1}\;v^S_M\,.
$$
\end{theor}
The last formula coinsides with Smirnov's formula for the singlet solutions
of the qKZ equation at level $0$ given in \cite{S1}.

\begin{remark}
F.~Smirnov used another basis $\,\{\omega_M\}\,$ in his construction of
solutions of the qKZ equation:
\beq
\label{omegaM}
\omega_M\,=\,
\prod_{\tsize{k\not\in M\atop m\in M}}{z_k-z_m\over z_k-z_m+\h}\;v^S_M\,.
\eeq
Set $\,{M'_{\rm ext}=\{n-\ell+1\lc n\}}\,$.
Then $\;\omega_{M'_{\rm ext}}\!=\,v_{M'_{\rm ext}}$.
The main property of this basis is given by Corollary \ref{symbasis}.
\end{remark}
\begin{remark}
In the difference case the space of periodic functions $\Fq$ plays the role of
the set of admissible contours $\Ch$ in the differential case and the subspace
$\Fcq$ is an analogue of the set of cycles $\Cc$. By the formula \r{total} and
Corollary \ref{77} we observe that for $n=2\ell$, which is the case of zero
weight, the integrand in the right hand side of the formula \r{solWSbas} for
solutions of the qKZ equation is a total difference. So that, for any
$\,W\in\Fcq^{\ot\ell}$ we have
$$
\psi_W(z_1\lc z_n)\,=\,0\,.
$$
Therefore, to produce nonzero solutions of the qKZ equation we inevitably have
to consider difference analogues of unclosed admissible contours in
the differential case.
\end{remark}

\Section{Solutions of the qKZ equation from the representation theory}
\label{777}
The aim of this section is to investigate the solutions of (\ref{qKZ0}) which
is obtained from the representation theory of the centrally extended Yangian
double ($\Aya$) \cite{Kh}, \cite{IK}.  It is a Hopf algebra associated with
the $R$-matrix $R^\pm(u)=\rho^\pm(u) R(u)$, where $\rho^\pm(u)$ are certain
scalar factors and $R(u)$ is given by \r{R-mat} \cite{Kh}.

The algebra $\Aya$ posseses two dimensional evaluation representation $V_z$
and the level one infinite dimensional representation ${\cal H}$. In \cite{Kh},
\cite{IK} the intertwining operators $\Phi(\z): \H\to \H\otimes V_{y}$ and
$\Psi(z):\H\to V_z\ot\H$ are constructed.  We define the components of
the intertwining operators by $\Phi(z) v = \sum_\nu\Phi_\nu(z)v\otimes v_\nu$,
$\Psi(z) v = \sum_\ep v_\ep\ot\Psi_\ep(z)v$, where $v\in \H$ and $v_\ep\in V$,
$\ep=\pm$.

Let us consider $(V^{\ot n})_\ell$ and $(V^{\ot n'})_{\ell'}$-valued functions:
\bea
\label{kz01}
\psi^{K}(\b_1\lc\b_n;\z_1\lc\z_{n'})
&=&\sum_{\scr M\subset\{1\lc n\}\atop\scr \# M=\ell}
\Omega_M^{K}(\b_1\lc\b_n;
\z_1\lc\z_{n'})\> v_M\ , \nn\\
\tilde\psi_M(\b_1\lc\b_n;\z_1\lc\z_{n'})
&=&\sum_{\scr K\subset\{1\lc n'\}\atop\scr \# K=\ell'}
\Omega_M^{K}(\b_1\lc\b_n;
\z_1\lc\z_{n'})\> v_K\ ,
\nn
\eea
where $\Omega_M^K(z;y)$ is a certain function proportional to the ratio of
traces
$$
{\tr_{\cal H}\left(e^{\g d} \Psi_{\ep_n}(\b_n) \ldots \Psi_{\ep_1}(\b_1)
\Phi_{\nu_1}(\z_1) \ldots \Phi_{\nu_{n'}}(\z_{n'}) \right)\over
\tr_{\cal H}\left(e^{\g d}\right)}
$$
where $\;M\,=\,\{\)j\ |\ \ep_j=-\,$, $j=1\lc n\}\subset\{1\lc n\}\,$,
$\;K\,=\,\{\)i\ |\ \nu_i=+\,$, $i=1\lc n'\}\subset\{1\lc n'\}\;$
and the proportionality coefficient does not depend on $M$ and $K$,
cf.\ \cite{KLP}. Notice that the sets $M$ and $K$ can be empty.
The trace of the composition of the intertwining operators was calculated
in \cite{KLP}, \cite{C}. The formula for $\Omega_M^{K}(\b;\z)$ is
\bea
&
\Omega_M^K(\b_1\lc\b_n;\z_1\lc\z_{n'})\,=\,\dl_{n-2\ell,n'-2\ell'}\,
\prod_{j=1}^n\ee^{i\pi(\ell-\ell')\b_j/p}
\prod_{i=1}^{n'}\ee^{i\pi(\ell'-\ell)\z_i/p}\times\kern-3em
\nn\\
&\quad {}\times\,\int_{\bar C}\cdots\int_{\bar C}d\v_1\cdots d\v_\ell
\int_{\tilde C}\cdots\int_{\tilde C}d\u_1\cdots d\u_{\ell'}\
P_M(\v;\b) P_K(\u;\z) \nn\\
&\qquad{}\times\,
\prod_{s=1}^\ell \prod_{j=1}^n
\left[\Ga{{\b_j-\v_s\over\g}}\Ga{{\v_s-\b_j-\h\over\g}}\right]
\prod_{r=1}^{\ell'}\prod_{i=1}^{n'}
\left[\Ga{{\u_r-\z_i\over\g}}\Ga{{\z_i-\u_r+\h\over\g}}\right]\nn\\
&\qquad\times
\prod_{1\le s'<s''\le\ell}\left[
{\sin(\pi({\v_{s'}-\v_{s''}})/{\g})\over\Gm{(({\v_{s'}-\v_{s''}-\h})/{\g})}
\Gm{(1+(\v_{s''}-\v_{s'}-\h)/{\g})}}\right]\nn\\
&\qquad \times\prod_{1\le r''<r'\le\ell'}\left[
{\sin(\pi({\u_{r'}-\u_{r''}})/{\g})\over\Gm{(({\u_{r'}-\u_{r''}+\h})/{\g})}
\Gm{(1+(\u_{r''}-\u_{r'}+\h)/{\g})}}\right]\nn\\
&\qquad{}\times
\ {\dis \prod_{r=1}^{\ell'}\prod_{j=1}^n\sin(\pi({\u_r-\b_j})/{\g})
\prod_{s=1}^\ell\prod_{i=1}^{n'}\sin(\pi({\v_s-\z_i-\h})/{\g})
\over \dis
\prod_{s=1}^\ell\prod_{r=1}^{\ell'}\sin(\pi({\u_r-\v_s})/{\g})
\sin(\pi({\u_r-\v_s+\h})/{\g})}\ ,
\label{general-trace}
\eea
In the last line of \r{general-trace} we correct a misprint made in \cite{KLP}.
Formula \r{general-trace} can also be obtained from the corresponding
formula for the quantum affine algebra in \cite{JM1} by taking the scaling
limit. Particular specializations of the formula \r{general-trace} are given
in \cite{L}, \cite{N2}.
\goodbreak

For $M=\{m_1<\ldots <m_\ell\}$ the polynomial $P_M(\v;\b)$ is defined
by the formula:
$$
P_M(\v;\b)\,=\,\prod_{a=1}^\ell
\Big(\prod_{j>m_a}(\v_a-\b_j)\prod_{j<m_a} (\v_a-\b_j-\h)\,\Big).
$$

The contour $\bar C$ and $\tilde C$ are specified as follows. The contour
$\bar C$ for the integration over $\v_a$, $a=1\lc\ell$ is a simple curve
separating the sets of the points $\Cup_{j=1}^n(\b_j+\h-p\ZZ_{\geq0})$,
$\Cup_{b=1}^{\ell'}(\u_b+\h-p\ZZ_{>0})$, $\Cup_{b=1}^{\ell'}(\u_b-p\ZZ_{>0})$
and $\Cup_{j=1}^n(\b_j+p\ZZ_{\geq0})$,
$\Cup_{b=1}^{\ell'}(\u_b+\h+p\ZZ_{\geq0})$,
$\Cup_{b=1}^{\ell'}(\u_b+p\ZZ_{\geq0})$, and going from $-\infty$ to $+\infty$.
More precisely the contour $\bar C$ admits a parametrization $\rho:\RR\to\CC$
such that $\rho(u)\to\pm\infty$ and $\Im\rho(u)$ has finite limits as
$u\to\pm\infty$.  Similarly the contour $\tilde C$ separates the sets
of the points $\Cup_{i=1}^{n'}(\z_i-p\ZZ_{\geq0})$,
$\Cup_{a=1}^{\ell}(\v_a-\h-p\ZZ_{\geq0})$,
$\Cup_{a=1}^{\ell}(\v_a-p\ZZ_{\geq0})$ and
$\Cup_{i=1}^{n'}(y_i+\h+p\ZZ_{\geq0})$, $\Cup_{a=1}^{\ell}(\v_a-\h+p\ZZ_{>0})$,
$\Cup_{a=1}^{\ell}(\v_a+p\ZZ_{>0})$, and going from $-\infty$ to $+\infty$.

Notice that for $\ell'=0$ the contour $\bar C$ coincides with the contour $C$
used in the definition of the hypergeometric integral \r{pairing}.

Due to the commutation relations of the intertwining operators and the cyclic
property of the trace, the function $\psi^K(\b_1\lc\b_n;\z_1\lc\z_{n'})$ solves
the qKZ equation \r{qKZ0} at level $-2+p/\h$ with respect to the variables
$\b_j$, $j=1\lc n$ for any set $K$ and any values of the parameters $\z_i$,
$i=1\lc n'$. On the other hand the function
$\tilde\psi_M(\b_1\lc\b_n;\z_1\lc\z_{n'})$ solves the qKZ equation at level
$-2-p/\h$ with respect to the variables $\z_i$, $i=1\lc n'$ for any set $M$ and
any values of the parameters $\b_j$, $j=1\lc n$, see \cite{JM1}, \cite{KLP}.

 From now on we assume that $p=2\h$ and we consider the solutions
$\psi^K(\b_1\lc\b_n;\z_1\lc\z_{n'})$ of the qKZ equation at level zero.

It follows from \r{general-trace} that
$\psi^{K}(\cd;\z_1\lc\z_{n'})\in\left(V^{\ot n}\right)_\ell\ $ with
$2\ell=n-n'+2\ell'$. In general $\,\psi^{K}(\b_1\lc\b_n;\z_1\lc\z_{n'})\,$ does
not satisfy the highest weight condition \r{sing}. Therefore we decompose it
into the isotypic components with respect to the $\sltwo$ action:
\beq
\psi^{K}(\cd;\z_1\lc\z_{n'})\,=\sum_{j=0}^{n}\, {(\Sigma^-)}^j
\>\psi^{[K,j]}(\cd;\z_1\lc\z_{n'})\ ,
\label{decomp}
\eeq
$$
\Sigma^+ \psi^{[K,j]}\>=\>0\,,\qqquad \psi^{[K,j]}\in(V^{\ot n})_{\ell-j}\,.
$$
Since the operators $\,K_j(z_1\lc z_n)\,$, cf.\ \r{qKZ0}, commute with
the $\frak{sl}_2$ action on $V^{\ot n}$, each component
$\psi^{[K,j]}(\cd;\z_1\lc\z_{n'})$ is a solution of the qKZ equation \r{qKZ0}.
Our general aim is to find functions
$$
W[K,j;\z_1\lc\z_{n'}]\in \widehat{\cal F}_q^{\ot(\ell-j)}
$$
such that
$$
\psi^{[K,j]}(\cd;\z_1\lc\z_{n'})  =\psi_{W[K,j;\z_1\lc\z_{n'}]}\ ,
$$
where $\psi_{W[K,j;\z_1\lc\z_{n'}]}$ is defined by \r{solutW}.
Notice that a function $W[K,j;\z_1\lc\z_{n'}]$ is not uniquely defined,
see Proposition \ref{5*} and the remark before it.
\Par
Up to now the formula of
$W[K,j;\z_1\lc\z_{n'}]$ for a general $K$ is not yet found.
Below we give simple examples.

\begin{example}1
Let $n'=n-2\ell\ge 0$, $K=\varnothing$. Then we have
\beq
\label{111}
\vp{\bigg|^|}
\psi^\varnothing (\cd;\z_1\lc\z_{n'})=\psi_{W[\varnothing,0;\z_1\lc\z_{n'}]}\ ,
\eeq
where
$$
W[\varnothing,0;\z_1\lc\z_{n'}]\)(t_1\lc t_\ell)\,=\,c_1\;
\prod_{a=1}^\ell{\dis \ee^{2\pi i(2a-1)\v_a/\g}
\prod _{k=1}^{n-2\ell}\left(1+\ee^{2\pi i(\v_a-\z_k)/\g}\right)
\over\dis \prod_{j=1}^{n}\left(1-\ee^{2\pi i(\v_a-\b_j)/\g}\right)}\ ,
$$
$$
\vp{\bigg|^|}c_1\,=\,(-1)^{\ell(\ell-1)/2}\>2^{\ell(n-\ell+1)}\>
(\pi i\)\g)^{\ell(2n-\ell+1)/2}\,.
$$
For $n'=0$, the empty product over $k$ equals $1$. In this case
$W[\varnothing,0]\in {\cal F}_q$ and $\psi_{W[\varnothing,0]}\equiv 0$
due to Lemmas \ref{4} and \ref{7**}.
\end{example}
\vsk>
\begin{example}2
Let $2\ell=n$, $n'=2$, $\ell'=1$, and $K=\{1\}$ or $\{2\}$. We have
\beq
\label{333}
\vp{\bigg|^|}\psi^K (\cd;\z_1,\z_2)\,=\,\psi_{W[K,0;\z_1,\z_2]}+ \Sigma^-
\left(\psi_{W[K,1;\z_1,\z_2]} \right)\ ,
\eeq
where
\beq
\label{33c}
W[K,0;\z_1,\z_2]\)(t_1\lc t_\ell)\,=\,\sigma_K(\z_2-\z_1-\h)\>
W[K,1;\z_1,\z_2]\)(t_1\lc t_{\ell-1})\>\tilde W(t_\ell)\,,
\eeq
$$
\vp{\bigg|^|}\sigma_{\{1\}}=1\,,\qqquad \sigma_{\{2\}}=-1\,,
$$
$$
W[K,1;\z_1,\z_2]\)(t_1\lc t_{\ell-1})\,=\,c_2\;\prod_{a=1}^{\ell-1}
{\dis \ee^{2\pi i(2a-1)t_a/\g}\prod_{k=1}^2 \sk{1+\ee^{2\pi i(t_a-\z_k)/\g}}
\over \dis \prod_{j=1}^n\sk{1-\ee^{2\pi i(t_a-\b_j)/\g}}}\ .
$$
$$
\vp{\bigg|^|_|}c_2\,=\,i(-1)^{(\ell-1)(\ell-2)/2}\>2^{\ell(\ell+1)}
\>\pi^{2(1-\ell)}(\pi i\)\g)^{(3\ell^2+3\ell+4)/2}\,,
$$
and $\tilde W$ is an arbitrary function from $\Fq$ such that
$\,2^\ell p\>\bigl(\tilde W(+\infty)-\tilde W(-\infty)\bigr)\,=\,1\,.$
For example we can take $\tilde W(t)=(2^\ell p)^{-1}
\prod_{j=1}^n\,(1-\ee^{2\pi i(t-z_j)/p})^{-1}$.
\end{example}

Let us briefly discuss the meaning of Examples~1 and~2 in the context of
integrable models of quantum field theory.
If we set $h=-\pi i$, $p=-2\pi i$, the model which corresponds to
our consideration in this paper is the $SU(2)$ invariant Thirring model (ITM).
Following the idea developed in \cite{DFJMN} the function 
$\psi^{K}(z_1\lc z_n;y_1\lc y_{n'})$ is equal modulo a scalar factor
to the $n$-particle form factor of the operator specified by
$\Phi_{\nu_1}(y_1)\cdots \Phi_{\nu_{n'}}(y_{n'})$.

In \cite{L} Lukyanov has introduced certain local operators. Some of them,
up to normalizations, are
\bea
\Lambda_m(\z)\,&\kern-18pt&{}=\;{1\over2}
(\Phi_{\ep_1}(y+\pi i)\Phi_{\ep_2}(\z)\>+\>
\Phi_{\ep_2}(y+\pi i)\Phi_{\ep_1}(\z)),
\label{current}
\\
T(\z)\,&\kern-18pt&{}=\;{1\over2}
(\Phi_{+}(y+\pi i)\partial_{\z}\Phi_{-}(\z)\>-\>
\Phi_{-}(y+\pi i)\partial_{\z}\Phi_{+}(\z)),
\label{energy}
\eea
where $m=0,\pm1$, $\ep_j=\pm$ and $m=(\ep_1+\ep_2)/2$.  We will show that
the form factors of \r{current} and \r{energy} include the form factors
obtained in \cite{S1}, \cite{S4}.
\begin{remark}
We change the signs of the second terms in the definitions of $\Lambda_0(\z)$
and $T(\z)$ compared with \cite{L} since we consider the $S$-matrix with a
nonsymmetric crossing symmetry matrix \cite{S1} as opposed to \cite{L}.
\end{remark}

We first present the formulae for form factors obtained in \cite{S1}, \cite{S4}
in terms of the functions $\psi_W$. Define the funtions
$\,W_\sigma\in\Fcq^{\ot(\ell-1)}$ and
$\,\tilde W_\sigma\in\Fcq^{\ot(\ell-1)}\ot\Fq\,$,
$\,\sigma=\pm\,$, by
\bea
&\kern-18pt& W_\sigma(t_1\lc t_{\ell-1})\,=\,
\prod_{a=1}^{\ell-1}\>\Bigl(\>\ee^{-(2a\mp1)\v_a}\,
\prod_{j=1}^{n}\>\bigl(1-\ee^{-(\v_a-\b_j)}\bigr)^{-1}\>\Bigr)\,,
\nn
\\[2pt]
&\kern-18pt& \tilde W_\sigma(t_1\lc t_\ell)\,=\,i\,2^{-\ell-1}\pi^{-1}\,
W_\sigma(t_1\lc t_{\ell-1})\,
\prod_{j=1}^{n}\>\bigl(1-\ee^{-(\v_\ell-\b_j)}\bigr)^{-1}.
\nn
\eea
Let $f^{\tau}_\sigma$, $\tau=\pm,3\,$, $\sigma=\pm\,$, be the form factors
of $SU(2)$ currents in the lightcone coordinates, cf.\ page 38 in \cite{S4},
and  $f_{\mu\nu}$, $\mu,\nu=0,1\,$, the form factors of the energy-momentum
tensor, cf.\ page 106 in \cite{S1}. Then, using the formula \r{det-ar-spin} and
Theorem \ref{7final}, we have
$$
f^-_\sigma(z_1\lc z_n)\,=\,{\rm const}\,\prod_{i<j}\>\zeta(z_i-z_j)\,
\exp\Bigl(\){\ell-1-\sigma\over2}{\tsize\sum\limits_{j=1}^n z_j}\>\Bigr)\,
\psi_{W_\sigma}(z_1\lc z_n)\,,
$$
\bea
f_{\mu\nu}(z_1\lc z_n)\,=\,{\rm const}\,\prod_{i<j}\>\zeta(z_i-z_j)\,
\sum_{j=1}^n\>(\ee^{z_j}-(-1)^\nu\)\ee^{-z_j})\>
\exp\Bigl(\){\ell-2\over2}{\tsize\sum\limits_{j=1}^n z_j}\>\Bigr)\,\times{}
&\kern-18pt&\nn
\\
{}\times\,\Bigl(\psi_{\tilde W_+}(z_1\lc z_n)-
(-1)^\mu\>\exp\)\bigl(\){\tsize\sum\limits_{j=1}^n z_j\)}\bigr)\>
\psi_{\tilde W_-}(z_1\lc z_n)\Bigr)
&\kern-18pt& \,,
\nn
\eea
where $2\ell=n$, $\,\zeta(z)$ is defined in \cite{S1}, page 107, and
the constants are independent of $z_1\lc z_n$. Formulae for the form factors
$f^3_\sigma$ and $f^{+}_\sigma$ can be obtained by applying the operator
$\Sigma^{-}$ once or twice to $f^{-}_\sigma$:
$\,f^3_\sigma=\Sigma^{-}f^{-}_\sigma$,
$\,f^{+}_\sigma=(\Sigma^{-})^2f^{-}_\sigma$.  Notice that the difference
equation satisfied by form factors of the $SU(2)$ ITM differs from
the qKZ equation \r{qKZ0} by the sign $(-1)^{n/2}$, cf.\ (110) in \cite{S1}.
This explains the appearance of the factors
$\,\exp\bigl(\){\ell-1-\sigma\over2}\sum\limits_{j=1}^n z_j\)\bigr)\,$ and
$\,\exp\bigl(\){\ell-2\over2}\tsize\sum\limits_{j=1}^n z_j\)\bigr)\,$ which,
in principle, are not $2\pi i$ periodic functions of $z_1\lc z_n$.

Example~1 for $n=2\ell+2$, $n'=2$ shows that the $n$-particle form factor of 
the operator $\Lambda_{-1}(\z)$ is proportional to
$\psi_{W[\varnothing,0;\z_1,\z_2]}$.
Notice that 
$$
W_\sigma(t_1\lc t_\ell)\,=\,\sigma^\ell
\lim{\z\to-\sigma \infty}\ee^{(\sigma-1)\)\ell\)\z}
c_1^{-1}W[\varnothing,0;\z+\pi i,\z](t_1\lc t_\ell)\,.
$$

Thus $f^{-}_\sigma$ is obtained from the form factor of $\,\Lambda_{-1}(\z)\,$
by the specialization of the value of $y$.  Similarly, $f^3_\sigma$ can be
obtained from the form factor of $\Lambda_0(\z)$, see Example~2. In the same
way, it is possible to show from Example~2 that $f_{\mu\nu}$, $\mu,\nu=1,2$,
can be obtained from the form factor of $T(\z)$ as suitable linear combinations
of the limits $\lim{\z\to\pm\infty}T(\z)$.

\subsection*{Acknowledgments}
S.~Pakuliak would like to thank Prof.\ T.~Miwa for the invitation to visit
RIMS in March\,--\,April 1997 where this work was started. The research of
S.~Pakuliak was supported in part by grants RFBR-97-01-01041 and by Award
No.~RM2-150 of the U.S.\ Civilian Research \& Development Foundation (CRDF)
for the Independent States of the Former Soviet Union.

\Appendix
\Section{}
\label{AF}

\begin{Proof}{Proposition~\ref{basis}}
Due to Lemmas~\ref{belong} and \ref{res}, both $\,\{\)w_M^{cl}\)\}_{\#M=\ell}$
\,and $\,\{\)\tilde w_M^{cl}\)\}_{\#M=\ell}$ \,are families of linear
independent functions from the space $\Fcl$, because
$\,\Res{}g_M^{cl}(\hat z_M)\ne 0\,$ under assumptions of the proposition.
Hence,
$$
\dim\Fcl\,\geq\,{n\choose\ell}\,,
$$
and it suffices to prove the opposite inequality:
$$
\dim\Fcl\,\leq\,{n\choose\ell}\,.
$$
Consider the space
$\;{\dsize\Gcl=\bigl\{f(t_1\lc t_\ell)\ |\ f(t_1\lc t_\ell)\!
\prod_{1\le a<b\le\ell}\!(t_a-t_b)\,\in\Fcl\)\bigr\}}\;$
\vv.15>
which is obviously iso\-mor\-phic to $\,\Fcl$. The definition of $\,\Fcl\,$
implies that $\,{f\in\Gcl}\,$ iff $\>f\>$ has the following properties:
\Par
\Item $f(t_1\lc t_\ell)$ is a symmetric rational function with at most simple
poles at the hyperplanes $\,{t_a=z_m}\,$, $\,a=1\lc\ell$, $\,m=1\lc n$,
\vsk.2>
\iitem $f(t_1\lc t_\ell)\to 0\,$ as $\,t_1\to\infty\,$ and $\,t_2\lc t_\ell\,$
are fixed,
\vsk.2>
\iitem
$\res{\,t_1=z_m\!}\!\bigl(\res{\,t_2=z_m\!}f(t_1\lc t_\ell)\bigr)=0\,$
\,for any $\,m=1\lc n\,$.
\Par
Let $\Xg$ be the set of sequences $(m_1\lc m_\ell)$ such that
$\,{m_1\lc m_\ell\in\{1\lc n\}}\,$ are pairwise distinct. For any permutation
$\,\sigma\,$ of $\,1\lc\ell\,$ and any $\,{\mg=(m_1\lc m_\ell)\in\Xg}\,$
set
$$
\mg^\sigma\,=\,(m_{\sigma_1}\lc m_{\sigma_\ell})\,.
$$

\begin{lemma}
\label{AXcl}
{\sl Let a function $f(t_1\lc t_\ell)$ have the properties\/
{\rm i)\,--\,\,iii)}. Then it has the form:}
$$
f(t_1\lc t_\ell)\,=\>\sum_{\mg\in\Xg}\,
{A_\mg\over(t_1-z_{m_1})\ldots (t_\ell-z_{m_\ell})}
$$
for suitable constants $\,\{A_\mg\}\,$ such that $\,A_\mg=A_{\mg^\sigma}\,$
for any $\,\mg\in\Xg\,$ and any permutation $\,\sigma$.
\end{lemma}
\begin{proof}
The lemma can be proved by the induction on $\,\ell\,$ using
the partial fraction expansion of a rational function of one variable,
cf.\ the proof of Lemma~\ref{AX}.
\end{proof}
Since Lemma~\ref{AXcl} implies that $\;\dsize\dim\Gcl\,\leq\,{n\choose\ell}\,$,
\;Proposition~\ref{basis} is proved.
\end{Proof}

\begin{Proof}{Proposition~\ref{3}}
Due to Lemmas~\ref{belong} and \ref{res}, both $\,\{\)w_M\)\}_{\#M=\ell}$
\,and $\,\{\)\tilde w_M\)\}_{\#M=\ell}$ \,are families of linear independent
functions from the space $\Fl$, because $\,\Res{}g_M(\hat z_M)\ne 0\,$ under
the assumptions of the proposition. Hence,
$$
\dim\Fl\,\geq\,{n\choose\ell}\,,
$$
and it suffices to prove the opposite inequality:
$$
\dim\Fl\,\leq\,{n\choose\ell}\,.
$$
Consider the space
$\;{\dsize\Gc=\bigl\{f(t_1\lc t_\ell)\ |\ f(t_1\lc t_\ell)\!
\prod_{1\le a<b\le\ell}\!(t_a-t_b)\,\in\Fl\)\bigr\}}\;$
\vv.15>
which is obviously iso\-mor\-phic to $\,\Fl$. The definition of $\,\Fl\,$
implies that $\,{f\in\Gc}\,$ iff $\>f\>$ has the following properties:
\Par
\atem $f(t_1\lc t_\ell)$ is a symmetric rational function with at most simple
poles at the hyperplanes $\,{t_a=z_m}\,$, $\,a=1\lc\ell$, $\,m=1\lc n$,
\vsk.2>
\bitem $f(t_1\lc t_\ell)\to 0\,$ as $\,t_1\to\infty\,$ and $\,t_2\lc t_\ell\,$
are fixed,
\vsk.2>
\bitem
$\res{\,t=z_m\!}\!f(t,t+\h,t_3\lc t_\ell)=0\,$ \,for any $\,m=1\lc n\,$.
\Par
For any $\,\mg\in\Xg$, \,set
$$
h_\mg(t_1\lc t_\ell)\,=\,\prod_{a=1}^\ell\,{1\over t_a-z_{m_a}\!}\;
\prod_{b=1}^{a-1}\,{t_a-z_{m_b}-\h\over t_a-z_{m_b}}\;.
$$

\begin{lemma}
\label{AX}
{\sl Let a function $f(t_1\lc t_\ell)$ have the properties\/
{\rm a)\,--\,\,c)}. Then it has the form:}
\beq
\label{fAh}
f(t_1\lc t_\ell)\,=\>\sum_{\mg\in\Xg}\,A_\mg\>h_\mg(t_1\lc t_\ell)
\eeq
for suitable constants $\,\{A_\mg\}\,$ such that
for any $\,\mg\in\Xg\,$ and any permutation $\,\sigma$
\beq
\label{Atau}
A_{\mg^\sigma}\,=\,A_\mg\!\prod_{\tsize{1\le a<b\le\ell\atop\sigma_a>\sigma_b}}
{z_{m_{\sigma_a}}-z_{m_{\sigma_b}}-\h\over
z_{m_{\sigma_a}}-z_{m_{\sigma_b}}+\h}\;.
\eeq
\end{lemma}
\begin{proof}
Consider the partial fraction expansion of $f(t_1\lc t_\ell)$ as a function
of $t_1$:
$$
f(t_1\lc t_\ell)\,=\>\sum_{m=1}^n\,{f_m(t_2\lc t_\ell)\over t_1-z_m}\;.
$$
The function $f_m(t_2\lc t_\ell)$ has the properties:
\Par
\textitem{${\rm a}'$)} $f_m(t_2\lc t_\ell)$ is a symmetric rational function
with at most simple poles at the hyperplanes $\,{t_a=z_j}\,$,
$\,a=2\lc\ell$, $\,j=1\lc n$,
\vsk.2>
\textitem{${\rm b}'$)} $f_m(t_2\lc t_\ell)\to 0\,$ as $\,t_2\to\infty\,$ and
$\,t_3\lc t_\ell\,$ are fixed,
\vsk.2>
\textitem{${\rm c}'$)}
$f_m(z_m+\h,t_3\lc t_\ell)=0\;$ and
$\,\res{\,t=z_l\!}\!f_m(t,t+\h,t_4\lc t_\ell)=0\,$ \,for any $\,l=1\lc n\,$.
\Par\nt
Hence,
\vv-.5>
$$
f_m(t_2\lc t_\ell)\,=\sum_{\tsize{l=1\atop\;l\ne m}}^n\,
{f_{lm}(t_3\lc t_\ell)\>(t_2-z_m-\h)\over(t_2-z_l)\>(t_2-z_m)}\;.
$$
The function $f_{lm}(t_3\lc t_\ell)$ has the properties:
\Par
\textitem{${\rm a}''$)} $f_{lm}(t_3\lc t_\ell)$ is a symmetric rational
function with at most simple poles at the hyperplanes $\,{t_a=z_j}\,$,
$\,a=3\lc\ell$, $\,j=1\lc n$,
\vsk.2>
\textitem{${\rm b}''$)} $f_{lm}(t_3\lc t_\ell)\to 0\,$ as $\,t_3\to\infty\,$
and $\,t_4\lc t_\ell\,$ are fixed,
\vsk.2>
\textitem{${\rm c}''$)}
$f_{lm}(z_l+\h,t_4\lc t_\ell)=0\,$, $\;f_{lm}(z_m+\h,t_4\lc t_\ell)=0\;$ and
$\,\res{\,t=z_k\!}\!f_{lm}(t,t+\h,t_5\lc t_\ell)=0\,$ \,for any $\,k=1\lc n\,$.
\Par\nt
Hence,
\vv-.5>
$$
f_{lm}(t_3\lc t_\ell)\,=\sum_{\tsize{k=1\atop\;k\ne l,m\!}}^n\,
{f_{klm}(t_4\lc t_\ell)\>(t_3-z_l-\h)\>(t_3-z_m-\h)\over
(t_3-z_k)\>(t_3-z_l)\>(t_3-z_m)}\;,
$$
etc. Finally we obtain the formula \r{fAh}.
\Par
For any $\,\mg\in\Xg$, \,let $\,\hat z_\mg\in\CC^\ell\,$ be the point defined
by $\,{\hat z_\mg=(z_{m_1}\lc z_{m_\ell})}\,$. It is clear that
\tenline
$\,\Res{}h_\lg(\hat z_\mg)=0\;$ for $\,\lg\ne\mg\,$ and
\beq
\label{Rhz}
\Res{}h_\mg(\hat z_\mg)\,=\prod_{1\le a<b\le\ell}
{z_{m_a}-z_{m_b}+\h\over z_{m_a}-z_{m_b}}\;.
\eeq
Since $f(t_1\lc t_\ell)$ is a symmetric function, we have
$\,\Res{}\]f(\hat z_\mg)=\Res{}f(\hat z_{\mg^\sigma})\,$ for any
$\,\mg\in\Xg\,$ and any permutation $\,\sigma$. Therefore,
$$
A_\mg\>\Res{}h_\mg(\hat z_\mg)\,=\,
A_{\mg^\sigma}\>\Res{}h_{\mg^\sigma}(\hat z_{\mg^\sigma})\,,
$$
which coincides with \r{Atau} because of \r{Rhz}. The lemma is proved.
\end{proof}

Since Lemma~\ref{AX} implies that $\;\dsize\dim\Gc\,\leq\,{n\choose\ell}\,$,
\;Proposition~\ref{3} is proved.
\end{Proof}

\Section{}
\label{A3*}

\begin{Proof}{Lemma~\ref{3*}}
We give the proof only for the formula \r{wMk}.
The proof of the formula \r{wMkcl} is similar.

Fix a subset $\,{M=\{m_2<\ldots< m_\ell\}\subset\{1\lc n\}}\,$. Say that
$\,{a\prec k}\,$ if $\,{m_a<k}\,$ and $\,{a\succ k}\,$ if $\,{m_a>k}\,$.
Let $f_1\lc f_n$ be the following functions:
$$
f_k(t_1\lc t_\ell)\,=\,\Bigl(1-\,{t_1-z_k-\h\over t_1-z_k}\,\Bigr)
\prod_{1\le l<k}\!{t_1-z_l-\h\over t_1-z_l}
\prod_{\tsize{a=2\atop a\prec k}}^\ell\,(t_1-t_a+\h)\,
\prod_{\tsize{a=2\atop a\succ k}}^\ell\,(t_1-t_a-\h)\,,
\qquad k\not\in M\,,
$$
\bea
f_{m_b}(t_1\lc t_\ell)\,&\kern-18pt&{}=\,\Bigl(t_1-t_b-\h-(t_1-t_b+\h)\,
{t_1-z_{m_b}-\h\over t_1-z_{m_b}}\,\Bigr)\,\times{}
\nn
\\[3pt]
&\kern-18pt&{}\>\times\,\prod_{1\le l<m_b}\!{t_1-z_l-\h\over t_1-z_l}
\prod_{2\le a<b}(t_1-t_a+\h)\prod_{b<a\le\ell}(t_1-t_a-\h)\,,
\qqquad b=2\lc\ell\,,
\nn
\eea
so that
$$
\prod_{a=2}^\ell\,(t_1-t_a-\h)\,-
\prod_{a=2}^\ell\,(t_1-t_a+\h)\,\prod_{j=1}^n{(t_1-z_j-\h)\over(t_1-z_j)}\;=\,
\sum_{k=1}^n\,f_k(t_1\lc t_\ell)\,.
$$
Set $\,{m_1=0}\,$ and $\,{m_{\ell+1}=n+1}\,$ until the end of the proof.
Let $\,{k\not\in M}\,$. Then there is a unique $\,a\,$, $\,1\le a\le\ell+1$,
such that $\,m_a<k<m_{a+1}\,$, and we get
$$
f_k(t_1\lc t_\ell)\>g_M(t_2\lc t_\ell)\,=\,\h\,(-1)^{a-1}\>
g_{M\cup\{k\}}(t_2\lc t_a,t_1,t_{a+1}\lc t_\ell)\,,
$$
$$
\Asym\bigl(f_k(t_1\lc t_\ell)\>g_M(t_2\lc t_\ell)\bigr)\,=\,
\h\,w_{M\cup\{k\}}(t_1\lc t_\ell)\,.
$$
On the other hand, for any $\,a=2\lc\ell\,$, the product
$\,f_{m_a}(t_1\lc t_\ell)\>g_M(t_2\lc t_\ell)\,$ is invariant with respect to
the permutation of $t_1$ and $t_a$. \,Hence,
$$
\Asym\bigl(f_{m_a}(t_1\lc t_\ell)\>g_M(t_2\lc t_\ell)\bigr)\,=\,0\,,
$$
which completes the proof.
\end{Proof}

\Section{}
\label{proofmain}
We give here the proof of Theorem~\ref{main} in order to make the paper
selfcontained.

\begin{Proof}{Theorem~\ref{main}}
We identify a subset $\,M\subset\{1\lc n\}\,$ with the sequence of signs
$\,(\ep_1\lc\ep_n)\,$ by the rule:
$$
M\,=\,\{\)i\ |\ \ep_i=-\>\}\,,
$$
cf.\ \r{vM}. Abusing notations we set
$\,{g_{\ep_1\lc\ep_n}\)=\)g_M}\,$ and $\,{w_{\ep_1\lc\ep_n}\,=\,w_M}\,$ if
the sequence $\,{(\ep_1\lc\ep_n)}\,$ corresponds to the subset $\,M$.
We will indicate explicitly that the functions $g_M$ and $w_M$ depend on
$z_1\lc z_n$, that is, for $\,{\#M=\ell}\,$ we will write
$\,g_M(t_1\lc t_\ell;z_1\lc z_n)\,$ and $\,w_M(t_1\lc t_\ell;z_1\lc z_n)$.

Notice that the integration contour $C$ in the definition of the hypergeometric
integral \r{pairing} obeys conditions which depends on $z_1\lc z_n$.
To indicate this we will write $C(z_1\lc z_n)$.
\Par
Let us introduce the twisted shift operators $Z_1\lc Z_n$ as follows
$$
Z_mf(t_1\lc t_\ell;z_1\lc z_n)\,=\,f(t_1\lc t_\ell;z_1\lc z_m+p\lc z_n)\,
\prod_{a=1}^\ell\,{t_a-z_m-p\over t_a-z_m-\h-p}\;,
$$
which means
\bea
&\kern-18.1pt&\nn
\\[-3\baselineskip]
\kern-2em f(t_1\lc t_\ell;z_1\lc z_m+p\lc z_n)\,
\prod_{a=1}^\ell\phi(t_a;z_1\lc z_m+p\lc z_n)\,=\!{} &\kern-18.1pt&
\nn
\\
{}=\,Z_mf(t_1\lc t_\ell;z_1\lc z_n)\,
\prod_{a=1}^\ell\phi(t_a;z_1\lc &\kern-18.1pt& z_n)\,.
\nn
\eea
We denote by $D_a$ the operator $D$ acting on the variable $\,t_a$,
cf.\ \r{tot-dif}. Let $\Fc'$ be the space of functions $f(t)$ such that
$(f-D\]f)\in\Fc^{(\ell+1)}$. Set
$\,\widetilde\Fc=\bigl(\Fc^{(\ell+1)}+\Fc'\bigr)^{\ot\ell}$. For any
$a=1\lc\ell$, $\,f\in\widetilde\Fc$ and $W\in\Fq^{\ot\ell}$ we have
$\,I(W,D_af)=0\,$, due to Lemmas~\ref{4} and \ref{48}.
We define the components of $R(z)$ by
$$
R(z)_{\ep,\ep}^{\ep,\ep}\>=\>1\,,\qqquad
R(z)_{-\ep,-\ep}^{\ep,\ep}\>=\>{z \over z+\h}\;,\qqquad
R(z)_{\ep,-\ep}^{-\ep,\ep}\>=\>{\h \over z+\h}\;,\qqquad \ep=\pm\,.
$$

It is easy to see that the function $\,\psi_W(z_1\lc z_n)\,$ defined
by \r{solutW} is a solution of the qKZ equation \r{qKZ0} if the following
relations hold:
\bea
\label{Rsym}
&\kern-18pt&
w_{\ep_1\lc\ep_{i+1},\ep_{i}\lc\ep_n}(\cdot;z_1\lc z_{i+1},z_i\lc z_n)\,=
\\[3pt]
&\kern-18pt& \]\]{}=\sum_{\ep'_i,\ep'_{i+1}=\pm}
R(z_i-z_{i+1})_{\ep_{i},\ep_{i+1}}^{\ep'_{i},\ep'_{i+1}}\,
w_{\ep_1\lc\ep'_i,\ep'_{i+1}\lc\ep_n}(\cdot;z_1\lc z_i,z_{i+1}\lc z_n)
\nn
\eea
for any $i=1\lc n$, \,and
\beq
\label{cyclic}
Z_1w_{\ep_1,\ep_2\lc\ep_n}(\cdot;z_1,z_2\lc z_n)\>-\>
w_{\ep_2\lc\ep_n,\ep_1}(\cdot;z_2\lc z_n,z_1)\,=\,
\sum_{a=1}^\ell\>D_a f^{(a)}_{\ep_1\lc\ep_n}
\eeq
for some functions $f^{(a)}_{\ep_1\lc\ep_n}\in\widetilde\Fc$, $\,a=1\lc\ell$.
Indeed, let $K_1(z_1\lc z_n)$ be the operator introduced in \r{qKZ0}.
Then relations \r{Rsym} and \r{cyclic} imply that
$$
Z_1w_{\ep_1\lc\ep_n}\>-\!\!\sum_{\,\ep'_1\lc \ep'_n=\pm\!}\!
(K_1)_{\ep_1\lc\ep_n}^{\ep'_1\lc \ep'_n}\,w_{\ep'_1\lc \ep'_n}\,=\,
\sum_{a=1}^\ell\>D_a f^{(a)}_{\ep_1\lc\ep_n}\,.
$$
Notice that
\bea
\kern-2em\int_{C(z_1+p,z_2\lc z_n)}\kern-.6em
Z_1w_{\ep_1\lc\ep_n}\>(t_1\lc t_\ell;z_1\lc z_n)\>
\phi(t_a;z_1\lc z_n)\>W(t_1\lc t_\ell;z_1\lc z_n)\>dt_a\;={}\!&\kern-18pt&
\nn
\\
{}=\!\int_{C(z_1\lc z_n)}\kern-.55em
Z_1w_{\ep_1\lc\ep_n}\>(t_1\lc t_\ell;z_1\lc z_n)\>
\phi(t_a;z_1\lc z_n)\>W(t_1\lc t_\ell;z_1\lc z_n)\>&\kern-18pt& dt_a
\nn
\eea
for any $a=1\lc\ell$ and $W\in\Fq^{\ot\ell}$, because the integrand has no
poles at $t_a=z_1-2\h$ and $t_a=z_1+5\h$. Therefore, the claim follows from
Lemmas~\ref{4} and \ref{48}. Equation \r{qKZ0} for $j>1$ can be proved
similarly.
\Par
We first prove the relation \r{Rsym}. If $(\ep_i,\ep_{i+1})\neq(-,-)$, 
hen we have
\bea
&\kern-18pt&
g_{\ep_1\lc\ep_{i+1},\ep_{i}\lc\ep_n}(\cdot;z_1\lc z_{i+1},z_i\lc z_n)\,={}
\nn
\\[3pt]
&\kern-18pt& \qquad{}=\sum_{\ep'_i,\ep'_{i+1}=\pm}
R(z_i-z_{i+1})_{\ep_{i},\ep_{i+1}}^{\ep^\prime_{i},\ep^\prime_{i+1}}\,
g_{\ep_1\lc\ep'_i,\ep'_{i+1}\lc\ep_n}(\cdot;z_1\lc z_i,z_{i+1}\lc z_n)\,,
\nn
\eea
which implies \r{Rsym}. Assume that $(\ep_i,\ep_{i+1})=(-,-)$. Let $\,a\,$ be
such that $\,m_a=i\,$ and $\,m_{a+1}=i+1$. Then we have
\bea
&\kern-36pt & \kern-13pt
g_{\ep_1\lc\ep_{i+1},\ep_{i}\lc\ep_n}(t_1\lc t_\ell;z_1\lc z_{i+1},z_i\lc z_n)
-g_{\ep_1\lc\ep_n}(t_1\lc t_\ell;z_1\lc z_n)\,={}
\nn
\\[3pt]
&\kern-36pt &
{}\!\!=\!\!\!\prod_{\tsize{b=1\atop b\ne a,a+1}}^\ell\!\!\!\!
\biggl(\>{1\over t_b-z_{m_b}}\prod_{1\leq l<m_b}\!\!
{t_b-z_l-\h\over t_b-z_l}\>\biggr)\!\!\!\!
\prod_{\tsize{1\le b<c\le\ell\atop (b,c)\neq(a,a+1)}}\!\!\!
(t_b-t_c-\h)\;\times{}
\nn
\\[4pt]
&\kern-36pt& {}\times\,\prod_{1\le l<i}\!
{(t_a-z_l-\h)\>(t_{a+1}-z_l-\h)\over(t_a-z_l)\>(t_{a+1}-z_l)}\;
{(z_{i+1}-z_i)\>(t_a-t_{a+1}-\h)\>(t_a-t_{a+1}+\h)\over
(t_a-z_i)\>(t_a-z_{i+1})\>(t_{a+1}-z_i)\>(t_{a+1}-z_{i+1})}
\nn
\eea
which is symmetric with respect to the permutation of $t_a$ and $t_{a+1}$.
Therefore
$$
w_{\ep_1\lc\ep_{i+1},\ep_{i}\lc\ep_n}(t_1\lc t_\ell;z_1\lc z_{i+1},z_i\lc z_n)
\,=\,w_{\ep_1\lc\ep_n}(t_1\lc t_\ell;z_1\lc z_n)
$$
in this case. Relation \r{Rsym} is proved.
\Par
To prove the relation \r{cyclic} we observe that for $\ep_1=+$ we have
$$
Z_1g_{\ep_1,\ep_2\lc\ep_n}(\cdot;z_1,z_2\lc z_n)\,=\,
g_{\ep_2\lc\ep_n,\ep_1}(\cdot;z_2\lc z_n,z_1)\,,
$$
and for $\ep_1=-$ we have
\bea
Z_1g_{\ep_1,\ep_2\lc\ep_n}(t_1\lc t_\ell;z_1,z_2\lc z_n)\>+\>(-1)^\ell
g_{\ep_2\lc\ep_n,\ep_1}(t_2\lc t_\ell,t_1;z_2\lc z_n,z_1)\,=\!{}
&\kern-18.1pt&\nn
\\[4pt]
{}=\,D_1
Z_1g_{\ep_1,\ep_2\lc\ep_n}(t_1\lc t_\ell;z_1,z_2\lc &\kern-18.1pt& z_n)\,.
\nn
\eea
It is clear that
$\,{Z_1g_{\ep_1,\ep_2\lc\ep_n}(\cdot;z_1,z_2\lc z_n)\in\widetilde\Fc}$.
\,Since $\,w_M=\Asym g_M\,$, we obtain
\bea
Z_1w_{\ep_1,\ep_2\lc\ep_n}(t_1\lc t_\ell;z_1,z_2\lc z_n)\>-\>
w_{\ep_2\lc\ep_n,\ep_1}(t_1\lc t_\ell;z_2\lc z_n,z_1)\,=\!{} &\kern-18.1pt&
\nn
\\
{}=\sum_{\sigma\in\Sb_\ell}^{\vp1}\sgn(\sigma)\>D_{\sigma_1}
f_{\ep_1\lc\ep_n}(t_{\sigma_1}\lc &\kern-18.1pt& t_{\sigma_\ell})\,,
\nn
\eea
for some functions $\,{f_{\ep_1\lc\ep_n}\in\widetilde\Fc}$,
which completes the proof.
\end{Proof}
\Par

As a corollary of the proof of Theorem \ref{main} given above we can describe
the properties of the bases $\,\{v^S_M\}\,$ and $\,\{\omega_M\}\,$.
If $M$ corresponds to $(\ep_1,\ldots,\ep_n)$ we set
$v^S_{\ep_1,\ldots,\ep_n}=v^S_M$ and $\omega_{\ep_1,\ldots,\ep_n}=\omega_M$.

\begin{corol}\label{symbasis}
The bases $\,\{v^S_M\}\,$ and $\,\{\omega_M\}\,$ satisfy 
the following equations:
\bea
\label{Rsymbasis}
&\kern-18pt&
u_{\ep_1\lc\ep_{i+1},\ep_{i}\lc\ep_n}(z_1\lc z_{i+1},z_i\lc z_n)\,={}
\\[3pt]
&\kern-18pt& \]\]{}=\,P_{i,i+1}\>R_{ii+1}(z_i-z_{i+1})\>
u_{\ep_1\lc\ep_i,\ep_{i+1}\lc\ep_n}(z_1\lc z_i,z_{i+1}\lc z_n)
\nn
\eea
for any $i=1\lc n$, where $u$ stands for either $v^S$ or $\omega$ and
$P_{i,i+1}$ is the permutation operator acting on the $i$-th and $(i+1)$-th
components.
\end{corol}

\begin{Proof}{Corollary \ref{symbasis}}
Set
$$
w(t_1 \lc t_\ell;z_1 \lc z_n)=\sum_M\, w_M(t_1 \lc t_\ell;z_1 \lc z_n)\>v_M\,.
$$
Then the equation \r{Rsym} is equivalent to
\beq
\label{sym1}
w(t_1 \lc t_\ell;z_1 \lc z_{i+1},z_i \lc z_n)\,=\,
P_{i,i+1}\>R_{i,i+1}(z_i-z_{i+1})\>w(t_1 \lc t_\ell;z_1 \lc z_n).
\eeq
Since $\;\tilde v_M(z_1\lc z_n)=\Res{}w(\hat z_M;z_1\lc z_n)\,$,
cf.\ \r{S-basis}, we have
\bea
\pm &\kern-18pt&
\tilde{v}_{\ep_1\lc\ep_{i+1},\ep_{i}\lc\ep_n}(z_1\lc z_{i+1},z_i\lc z_n)
\,={}
\nn
\\[3pt]
&\kern-18pt&\!\!{}=\,P_{i,i+1}\>R_{i,i+1}(z_i-z_{i+1})\>
\tilde{v}_{\ep_1\lc\ep_{i},\ep_{i+1}\lc\ep_n}(z_1\lc z_{i},z_{i+1}\lc z_n),
\nn
\eea
with the sign $-$ for the case $(\ep_i,\ep_{i+1})=(-,-)$ and $+$, otherwise.
Corollary \ref{symbasis} follows from this relation and formulae \r{S-basis2},
\r{omegaM}.
\end{Proof}

One can show that the bases $\,\{v^S_M\}\,$ and $\,\{\omega_M\}\,$ are uniquely
determined by the equation \r{Rsymbasis} and the respective normalizations
$\;{v^S_{M_{\rm ext}}\!=\,v_{M_{\rm ext}}^{\vp1}\!+\ldots{}}\,$,
$\ \omega_{M'_{\rm ext}}\!=\,v_{M'_{\rm ext}}\,$.

\Section{}
\label{DQ}
\begin{Proof}{Lemma~\ref{7procent}}
Fix a subset $\,M\subset\{1\lc n\}$, $\,\#M=\ell\>$. Consider a function
\bea
\label{fty}
\llap{$f(t,y)\,$}&\kern-18pt&{}=\;{P^+_M(t)\>P^-_M(t)\over t-y}\;-\;
{P^+_M(t+\h)\>P^-_M(t+\h)\over t-y+2\h}\ -
\\[3pt]
&\kern-18pt&{}-\;{\h\>P^+_M(y)\>P^-_M(t)\over(t-y)\>(t-y+\h)}\;-\;
{\h\>P^+_M(t+\h)\>P^-_M(y-\h)\over(t-y+\h)\>(t-y+2\h)}\;.
\nn
\eea
It is easy to see that $\,f(t,y)\,$ is a polynomial in $t,y$.
Moreover, for any $\,m\in M\,$ we have
\beq
\label{fzm}
f(t,z_m+2\h)\,=\,D\Bigl(\prod_{\tsize{k=1\atop k\neq m}}^n\,(t-z_k-2\h)\Bigr)
\,-\,\h\prod_{\tsize{k=1\atop k\neq m}}^n\,(z_m-z_k-\h)\>\mu_M^{(m)}(t)\,.
\eeq
Let
\beq
\label{fgPq}
g(t,y)\,=\,\Bigl[{f(t,y)\over P^+_M(y)}\Bigr]_{+,y}\,,\qqquad
q(t,y)\,=\,f(t,y)-P^+_M(y)\>g(t,y)
\eeq
where we take the polynomial part with respect to $y$. Then $\,q(t,y)\,$ is
a polynomial of degree less than $\ell$ with respect to $y$. We define
polynomials $\,q^{(1)}(t)\lc q^{(\ell)}(t)\,$ by the rule:
\beq
\label{qa}
q(t,y)\,=\,\sum_{a=1}^\ell\,q^{(a)}(t)\,y^{a-1}\,.
\eeq
Since $\,P^+_M(z_m+2\h)=0\,$, using formulae \r{fzm} and \r{fgPq} we obtain
$$
\kern-10pt
D\Bigl(\prod_{\tsize{k=1\atop k\neq m}}^n\,(t-z_k-2\h)\Bigr)\,=\,
\h\,\prod_{\tsize{k=1\atop k\neq m}}^n\,(z_m-z_k-\h)\>\mu_M^{(m)}(t)\,+\,
\sum_{a=1}^{\ell}\,q^{(a)}(t)\>(z_m+2\h)^{a-1}\,,
$$
Rewrite $f(t,y)$ in the form:
$$
f(t,y)\,=\,P^+_M(t+\h)\bigl(h(t,y)-h(t+\h,y)\bigr)\,+\,
P^-_M(t)\bigl(\tilde h(t,y)-\tilde h(t+\h,y)\bigr)\,,
$$
where
$$
h(t,y)\,=\;{P^-_M(t)-P^-_M(y-\h)\over t-y+\h}\;,\qqquad
\tilde h(t,y)\,=\;{P^+_M(t)-P^+_M(y)\over t-y}\;.
$$
Let
$$
q(t,y)\,=\,P^+_M(t+\h)\bigl(r(t,y)-r(t+\h,y)\bigr)\,+\,
P^-_M(t)\bigl(\tilde r(t,y)-\tilde r(t+\h,y)\bigr)
$$
be the corresponding decomposition of $q(t,y)$, cf.\ \r{fgPq}, where
$$
r(t,y)\,=\,h(t,y)\,-\,P^+_M(y)\>\Bigl[{h(t,y) \over P^+_M(y)}\Bigr]_{+,y}\;,
\qqquad\tilde r(t,y)\,=\;\tilde h(t,y)\,-\,P^+_M(y)\>
\Bigl[{\tilde h(t,y)\over P^+_M(y)}\Bigr]_{+,y}\,.
$$
Obviously we have
$$
\tilde r(t,y)\,=\;{P^+_M(t)-P^+_M(y)\over t-y}\;=\;
\Bigl[{P^+_M(t)\over t-y}\Bigr]_{+,t}\;.
$$
To rewrite appropriately $\,r(t,y)\,$ we use Lemma~\ref{t-y} below for
replacing a polynomial part with respect to $y$ by a polynomial part
with respect to $t$:
\bea
r(t,y)\,&\kern-18pt&{}=\;{P^-_M(t)-P^-_M(y-\h)\over t-y+\h}\;+\;
P^+_M(y)\>\Bigl[{P^-_M(y-\h)\over P^+_M(y)\>(t-y+\h)}\Bigr]_{+,y}\,={}
\nn
\\[3pt]
&\kern-18pt&{}=\,\Bigl[{P^-_M(t)\over t-y+\h}\Bigr]_{+,t}-\;
P^+_M(y)\>\Bigl[{P^-_M(t)\over P^+_M(t+\h)\>(t-y+\h)}\Bigr]_{+,t}\,={}
\nn
\\[3pt]
&\kern-18pt&{}=\,\Bigl[{P^-_M(t)\over P^+_M(t+\h)}\cdot
{P^+_M(t+\h)-P^+_M(y)\over t-y+\h}\Bigr]_{+,t}\,=\,\Bigl[{P^-_M(t)\over
P^+_M(t+\h)}\>\Bigl[{P^+_M(t+\h)\over t-y+\h}\Bigr]_{+,t}\>\Bigr]_{+,t}\,.
\nn
\eea
Finally
$$
q(t,y)\,=\,P^+_M(t+h)\>\Bigl[\>T_\h\Bigl({P^-_M(t)\over P^+_M(t+h)}\,
\Bigl[{P^+_M(t+h)\over t-y+\h}\Bigr]_{+,t}\Bigr)\Bigr]_{+,t}\>+\;
P^-_M(t)\>\Bigl[\>T_\h\Bigl({P^+_M(t)\over t-y}\Bigr)\Bigr]_{+,t}\,.
$$
Therefore, the polynomials $\,q^{(1)}\lc q^{(\ell)}\,$ defined by \r{qa}
coincide with the polynomilas $\,Q_M^{(1)}\lc Q_M^{(\ell)}\,$ given by \r{QM}.
Lemma~\ref{7procent} is proved.
\end{Proof}

\begin{lemma}
\label{t-y}
For any rational function $\,f(u)\,$ we have
$$
\Bigl[{f(u)\over u-x}\Bigr]_{+,u}=\,\Bigl[{f(x)\over x-u}\Bigr]_{+,x}\,.
$$
\end{lemma}
\vsk.5>
\begin{proof}
\hfill$\dsize
\Bigl[{f(u)\over u-x}\Bigr]_{+,u}=\,\biggl[{[f]_+(u)\over u-x}\biggr]_{+,u}=
\;{[f]_+(u)-[f]_+(x)\over u-x}\;=\,
\biggl[{[f]_+(x)\over x-u}\biggr]_{+,x}=\,\Bigl[{f(x)\over x-u}\Bigr]_{+,x}\,$.
\end{proof}

\end{document}